\documentclass[11pt,oneside,reqno]{amsart}
\usepackage{amsmath}
\usepackage{amsfonts}
\usepackage{amssymb}
\usepackage{amscd}
\usepackage{mathrsfs}

\catcode `\@=11
\def\numberbysection{\@addtoreset{equation}{section}
         \renewcommand{\theequation}{\thesection.\arabic{equation}}}
\numberbysection
\def\subsubsection{\@startsection{subsubsection}{3}%
  \normalparindent{.5\linespacing\@plus.7\linespacing}{-.5em}%
  {\normalfont\bfseries}}
\def\iso{\cong}
\def\a{\alpha}
\def\b{\beta}

\def\eps{\epsilon}
\def\l{\lambda}

\def\nn{\nonumber}
\def\la{\langle}
\def\ra{\rangle}
\def\Z2{\mathbb {Z}/2\mathbb{Z}}

\def\Hom{\mathrm{Hom}}
\def\C{\mathscr{C}}
\def\F{\mathscr{F}}
\def\Forb{\F_{\rm stringy}}
\def\calF{{\mathcal F}}

\def\dkg{D(k[G])}
\def\dbkg{D^{\beta}(k[G])}
\def\subdbkg{\dbkg^{\rm comm}}
\def\subdkg{\dkg^{\rm comm}}

\def\kg{\field[G]}

\def\isoto{\stackrel{\sim}{\rightarrow}}

\def\e{{\rm Eu}}

\def\etop{\e}

\def\hotimes{\widehat{\otimes}}
\def\hhotimes{\widehat{\hotimes}}
\def\Xm{X^m}

\def\bm{{\bf m}}
\def\F{{\mathcal F}}
\def\R{{\mathcal R}}

\newcommand{\eu}[1]{\e_{#1}}
\newcommand{\eut}[1]{\e_{#1,t}}

\def\L{\mathscr L}
\def\G{\mathscr G}

\def\D{\dkg-Mod} 

\def\unit{\upsilon}
\def\unitc{{\mathbb I}_{\mathcal C}}
\def\twist{{\mathbb I}_{\chi}}
\def\counit{\epsilon}

\def\Kglobal{\Korb}
\def\Korb{K_{\rm global}}
\def\Ksmall{K_{\rm small}}
\def\Kfull{K_{\rm full}}
\def\dprind{\rm Ind^{DPR}}
\def\ind{\rm Ind}
\def\res{\rm Res}
\def\bm{{\mathbf m}}
\def\Ag{{}^Ag}
\def\Bg{{}^Bg}
\def\Cg{{}^Cg}
\def\piAa{\pi^A_{\alpha}}

\def\Ch{\mathfrak{C}h}

\def\field{k}

\newcommand{\dkgv}[1]{D^{#1}(k[G])}

\newcommand{\elt}[2]{\, #1 \underset{\displaystyle #2}{\llcorner}}
\newcommand{\eltk}[2]{1_{#1,#2}}

\def\gx{\elt{g}{x}}
\def\hy{\elt{h}{y}}

\newtheorem{thm}{Theorem}[section]
\newtheorem{lem}[thm]{Lemma}
\newtheorem{prop}[thm]{Proposition}
\newtheorem{cor}[thm]{Corollary}

\newtheorem{dfprop}[thm]{Definition-Proposition}

\theoremstyle{definition}

\newtheorem{df}[thm]{Definition}

\newtheorem{rmk}[thm]{Remark}
\newtheorem{nota}[thm]{Notation}
\newtheorem{ex}[thm]{Example}
\newtheorem*{assumption}{Assumption}
\newtheorem*{cav}{Caveat}

\newcommand{\cf}{{\mathscr F}}
\newcommand{\ccx}{{\ci\cx}} 
\newcommand{\ci}{{\mathfrak I}}
\newcommand{\cx}{\mathfrak{X}}

\def\uk{k^*}

\def\morita{\sim_{\rm Morita}}
\begin{document}

\title{The  Drinfel'd Double and Twisting in Stringy Orbifold Theory}

\author[Ralph M. Kaufmann]{Ralph M. Kaufmann\\
Purdue University\\}

\address{Ralph Kaufmann:
 Department of Mathematics, Purdue University
150 N. University Street, West Lafayette, IN 47907-2067}
\email{rkaufman@math.purdue.edu}
\author[David Pham]{David Pham\\
University of Connecticut}
\address{David Pham:  Department of Mathematics, University of Connecticut,
196 Auditorium Rd, Storrs, CT 06269-3009}
\email{pham@math.uconn.edu}
\begin{abstract}
This paper exposes the fundamental role that the Drinfel'd double $\dkg$
of the group ring of a finite group $G$ and its
twists $\dbkg$, $\beta \in Z^3(G,\uk)$
as defined by Dijkgraaf--Pasquier--Roche play in stringy orbifold theories
 and their twistings.

The results pertain to three different aspects of the theory.
First, we show that $G$--Frobenius algebras  arising in
global orbifold cohomology or K-theory  are most naturally
 defined as elements in the braided category of $\dkg$--modules.
Secondly, we obtain a geometric realization of the Drinfel'd double
as the global orbifold $K$--theory of global quotient given by the
inertia variety of a point with a $G$ action on the one hand and
more stunningly a geometric realization of its representation ring
in the braided category sense as the full $K$--theory of the stack
$[pt/G]$. Finally, we show how one can use the co-cycles $\beta$
above to twist a) the global orbifold $K$--theory of the inertia of
a global quotient and more importantly b) the stacky $K$--theory of
a global quotient $[X/G]$. This corresponds to twistings with a
special type of 2--gerbe.
\end{abstract}
\maketitle


\subsection*{Introduction}
The Drinfel'd double $\dkg$ of a group ring of a finite group $G$ and in particular its twisted version $\dbkg$ where  $\beta \in Z^3(G,\uk)$ were introduced
and studied by Dijkgraaf, Pasquier and Roche \cite{DPR} (see \cite{ACM} for
a very nice brief summary).
Their aim was to understand the constructions of \cite{DVVV} concerning
 orbifold conformal field theory on the one hand and the constructions of \cite{DW} pertaining to orbifold Chern--Simons theory on the other. We will
realize these algebraic constructions geometrically using the
orbifold $K$--theory of \cite{JKK2} and newly defined twists.

Mathematically, the appearance of the Drinfel'd double $\dkg$ as a
main character in orbifold theory has its roots in \cite{FQ,F} where
a 2+1 dimensional theory was considered. See also \cite{FHT1,FHT2}
for related material on equivariant $K$--theory of a compact group
$G$. The importance and algebraic relevance of $\dkg$ in the theory
of  $G$--Frobenius algebras was made precise in \cite{disc} where we
showed that any $G$--Frobenius algebra is a $\dkg$--module and in
particular also a $k[G]$--module algebra and $k[G]$ co--module
algebra. $G$--Frobenius algebras arise in the 1+1 dimensional theory
(\cite{orb}) such as orbifold Gromov--Witten theory \cite{CR} and
hence in orbifold cohomology \cite{CR,FG} in particular. In this
paper, we go one step further and give a definition of a
$G$--Frobenius algebra and more generally a $G$--Frobenius object in
terms of the braided tensor category $\D$ of $\dkg$--modules. The
rather lengthy original definition of a $G$--Frobenius algebra
\cite{wisc,orb} then can be replaced by  the statement that a
$G$--Frobenius object is a Frobenius object in $\D$ which satisfies
two additional axioms (S) and (T) of which the former is the famous
trace axiom. This is the content of Theorem \ref{frobthm}.

Another upshot of the categorical
treatment is that these objects give the right algebraic structure
to encode the trace axiom in infinite dimensional situations.
We recall that in \cite{JKK2}, we introduced
pre--Frobenius algebras with trace elements to be able to write the trace
axiom. This was necessary since
the Chow ring of a smooth projective variety
 may not be a Frobenius algebra as it can be
infinite dimensional.
Here by a Frobenius algebra we
 mean a unital associative commutative algebra with a non--degenerate
even symmetric invariant bi--linear pairing. Nevertheless there are
traces one can define using the trace elements and for these the
trace axiom holds. In the categorical context any Frobenius object
defines a trace for any endomorphism which we call Frobenius trace
or F--trace for short. In particular the trace elements of
\cite{JKK2} can be recovered as the F--traces of the relevant
endomorphisms. This fact holds true in all the known constructions
involving the string and global versions of the  functors $\F\in
\{H^*,K^*,A^*,K_0\}$ \cite{FG,AGV,CR,JKK2} which is shown in Theorem
\ref{frobobthm} Thus $\dkg$ is at the bottom of the very definition
of the algebras associated to global orbifolds. Analogous statements
are true for singularities with symmetries \cite{wisc,orb,singorb}.

The Drinfel'd double makes its appearance in two more guises. First
we show that in the case of an Abelian
symmetry group $G$ the global $K$--theory as defined in \cite{JKK2}, see also Section 3 for a review, of the inertia variety of a point with the trivial $G$ action satisfies $\Kglobal^*(I(pt,G),G)=\dkg$ as an algebra.
In the non--Abelian case the resulting algebra together with its $G$--action
is Morita equivalent to $\dkg$ as a groupoid, see Corollary \ref{moritacor}.

The most stunning appearance of $\dkg$ is the one of Theorem
\ref{ptdprthm} were we prove that $\Kfull([pt/G])\iso Rep(\dkg)$.
Here the non--commutativity in the ring structure is now given by
the natural braiding of the moniodal category of representations.

Armed with these results, we define twistings
by co--cycles in $Z^i(G,\uk)$ where $i=1,2,3$ for
the various theories associated to a global quotient $(X,G)$. That is
in other words twist by 0--,1--, and 2--gerbes that are pulled
back from $[pt/G]$ or gerbes on $X$ that are trivial but not equivariantly
trivial see \cite{thaddeus} and \cite{hitch} for this point of view
of gerbes.

The $0$ twists are performed on $\Kglobal(X,G)$ or any of the other stringy
functors $\F$. They correspond to the Ramond twist defined in \cite{wisc,orb}.
The twists by $1$--gerbes are identified as the twist of
discrete torsion that were algebraically defined in \cite{disc}.
Finally the most interesting twists come from $2$--gerbes.
There are basically two types. First we can transgress the $2$--gerbe to the
inertia variety $I(X,G)$ considered together with its $G$--action and
 then consider twists on $\Kglobal((I(X,G),G))$. Here the
twist will just be a special type of discrete torsion.
 However, we do recover the
algebra structure of $\dbkg$ for the $\b$ twisted $K^{\beta}((I(pt,G),G)$.
The more intriguing twist is on $\Kfull[X/G]$.
We would like to note that in \cite{AR2}  a different interesting
 twist on $\Kfull(\ccx)$, that is the orbifold $K$--theory of the inerita
orbifold $\ccx$ for an orbifold $\cx$, was considered. In our case
we remain on $\Kfull(\cx)$ and our twist yields the natural
generalization of the results above. Namely $K^{\beta}([pt/G])\iso
Rep(\dbkg)$ see Theorem \ref{dbkgrepthm}. This result is striking in
several aspects. The most prominent feature being that the
representation ring of $\dbkg$ is understood in the braided monoidal
setting with a non--trivial associator.
 This tells us that this twist twists outside
the associative world. {\em A posteriori} this is however not
totally  unexpected, since we know from the work of Moore and
Seiberg \cite{MooreSei} that the fusion ring is not associative in
general, but only associative in the braided monoidal category
sense. We can of course get an associative algebra by restricting to
the dimensions of the intertwiners and defining a Verlinde algebra,
see Section 4 and also \cite{FHT1,FHT2} for related material.

The paper is organized as follows:

In Section 1, we review all the necessary definitions for the
twisted Drinfel'd double including DPR induction and the relevant
background from braided monoidal categories. Section 2  contains the
first set of results that pertain to the definition of
$G$--Frobenius algebra objects. The third section starts with a
brief review of the constructions of \cite{JKK2} and introduces all
the variants of stringy $K$--theory we will consider. Section 3
terminates with the second and third appearance of the Drinfel'd
double:  a) as the global $K$--theory of the inertia of $(pt,G)$ and
b) in the Theorem that $\Kfull([pt/G])\iso Rep(\dkg)$. The various
twistings are  contained in Section 4. Here we consider twists of
0--,1-- and 2-- gerbes on global quotients that are trivial but not
equivariantly trivial.

\section*{Acknowledgments}
It is a pleasure for R.K. to thank the Mittag--Leffler Institute and
the Max--Planck--Institute for Mathematics for their hospitality. We
would also like to thank Takashi Kimura,  Tyler Jarvis and Michael
Thaddeus for very useful discussions.

\section{The twisted Drinfel'd double}
In this section, we collect the basic definitions and constructions
of the twisted Drinfel'd double for the readers' convenience.

\subsection{Basic definitions}
\begin{df} For a finite group $G$ and an element $\b\in
Z^3(G,\uk)$, the {\it twisted Drinfel'd double} $\dbkg$ is the
quasi-triangular quasi-Hopf algebra whose

\begin{itemize}
\item[(1)] underlying {\em vector space has} the basis $\elt{g}{x}$ with $x,g\in
G$
\begin{equation}
\dbkg=\bigoplus k \, \gx
\end{equation}
\item[(2)] algebra structure is given by
\begin{equation}
\gx\hy=\delta_{g,xhx^{-1}}\theta_{g}(x,y)\elt{g}{xy}
\end{equation}
 where
 \begin{equation}
\label{thetaeq}
\theta_{g}(x,y)=\frac{\beta(g,x,y)\beta(x,y,(xy)^{-1}g(xy))}{\beta(x,x^{-1}gx,y)}
\end{equation}

\item[(3)] co-algebra structure is given by
\begin{equation}
\Delta(\gx)=\sum_{g_1 g_2=g} \gamma_x(g_1,g_2) \elt{g_1}{x}\otimes
\elt{g_2}{x}
\end{equation}
where

\begin{equation}
\label{gammadefeq}
\gamma_x(g_1,g_2)=\frac{\beta(g_1,g_2,x)\beta(x,x^{-1}g_1x,x^{-1}g_2x)}{\beta(g_1,x,x^{-1}g_2x)}.
\end{equation}

\item[(4)] The Drinfel'd associator $\Phi$ is given by
\begin{equation}
\Phi=\sum_{g,h,k\in G}\beta(g,h,k)^{-1} \elt{g}{e}\otimes \elt{h}{e}
\otimes \elt{k}{e}
\end{equation}
\item[(5)]
The $R$ matrix is given by

\begin{equation}
R=\sum_{g\in G} \elt{g}{e}\otimes \elt{\mathbf{1}}{g}, \text{ where
} {\elt{\mathbf{1}}{g}}=\sum_{h\in G}\elt{h}{g}
\end{equation}

\item[(5)]
The antipode $S$ is given by
\begin{equation}
S(\gx)=\frac{1}{\theta_{g^{-1}}(x,x^{-1})\gamma_x(g,g^{-1})} \text {
} \elt{x^{-1}g^{-1}x}{x^{-1}}
\end{equation}

\end{itemize}
\end{df}

\begin{rmk} There are several things which we would like to point
out:
\begin{enumerate}
\item In case $\beta\equiv 1$, that is $\beta$ is trivial, we
obtain the a braided Hopf algebra  $\dkg$ which is the Drinfel'd
double of the group ring.
\item The algebra is associative and the unit of this algebra is $\elt{\mathbf 1}{e}$

\item There is an injection of algebras $k[G]^*\to \dbkg$ given by $\delta_g\mapsto
\elt{g}{e}$, where $\delta_g(h):=\delta_{g,h}$, since
\begin{equation} \label{projectors}
\elt{g}{e}\elt{h}{e}=\delta_{g,h}\elt{g}{e}
\end{equation}

\item There is a special element $v^{-1}$ which is central.
It is given by
\begin{equation}
\label{veq}
v^{-1}=\sum_{g\in G}\elt{g}{g}
\end{equation}
In case $\beta\equiv 1$ this is the element which gives the inner
operation of $S^2$ of the braided Hopf-algebra $\dkg$ \cite{Kassel}.
\item the various $\theta_g$ are almost co-cycles for $G$
\begin{equation}
\theta_g(x,y)\theta_g(xy,z)= \theta_g(xy,z)\theta_{x^{-1}gx}(y,z)
\end{equation}
it follows that when $\theta_g$ is restricted to $Z(g)\times Z(g)$
it  {\em is} a 2--co-cycle for $Z(g)$.
\end{enumerate}
\end{rmk}

\subsection{The braided monoidal category $\dbkg-Mod$}
\label{braidedpar}
Since $\dbkg$ is a quasi-triangular quasi Hopf algebra, there is a
natural braided monoidal structure on the category of its modules.
We recall that if $U$ and $V$ are modules over $\dbkg$ or in general
any quasi--triangular quasi--Hopf algebra
$(H,\mu,\eta,\Delta,\eps,S,\Phi,R)$ then $U\otimes V$ has the
structure of an $H$ module via $\Delta: H \to H\otimes H$.

Recall (see e.g.\ \cite{Kassel}) that in general for three
representation $U,V,W$ and elements $u\in U, v\in V, w\in W$ the
associator is given by

\begin{eqnarray}
a_{U,V,W}: (U\otimes V) \otimes W &\to &U\otimes (V\otimes W)\nn\\
 a_{U,V,W}((u\otimes v)
\otimes w)&=&\Phi(u\otimes(v\otimes w))
\end{eqnarray}
and likewise for two representations $U,V$ and elements $u\in U,
v\in V$ the braiding is given by
\begin{eqnarray}
c_{U,V}: U\otimes V&\to &V\otimes U\nn\\
 c_{U,V}(u\otimes v)&=&\tau_{U,V}(R(u\otimes v))
\end{eqnarray}
where $\tau_{U,V}(u\otimes v)=v\otimes u$.

In particular, let $U,V,W$ be $\dbkg$ modules and let $u_g\in U_g$,
$v_h\in V_h$, $w_k\in W_k$ be homogeneous elements with respect to
the grading by $G$ then

\begin{equation}
\label{associate} a_{U,V,W}((u_g\otimes v_h) \otimes
w_k)=\beta^{-1}(g,h,k) u_g\otimes (v_h  \otimes w_k)
\end{equation}
and
\begin{equation}
\label{braid}
 c_{U,V}(u_g\otimes v_h)=\rho(\elt{ghg^{-1}}{g})(v_h)\otimes
 u_g=\phi(g)(v_h)\otimes u_g
\end{equation}

Moreover on $U\otimes V$ the $\dbkg$ module structure is given by
\begin{equation}
\label{action} \rho(\gx)(u_h\otimes
v_k)=\delta_{xhkx^{-1},g}\gamma_x(xhx^{-1},xkx^{-1})
\rho(\elt{xhx^{-1}}{x})(u_h)\otimes \rho(\elt{xkx^{-1}}{v_k})
\end{equation}

\begin{rmk}
It is well known \cite{MooreSei} that the pentagon relation for
associativity constraint is equivalent to the fact that $\beta$ as a
function on $G^3$ is an element of $Z^3(G,\uk)$.
\end{rmk}

\begin{prop}
\label{gradingprop} Any left $\dbkg$-module $(\rho,A)$ is $G$ graded
$A=\bigoplus_{g\in G}A_g$ and if $\pi_g$ denotes the projection of
$A$ onto $A_g$ then
\begin{enumerate}
\item $\rho(\elt{g}{e})=\pi_g$
\item $\rho(\gx)= \pi_g \circ \rho(\gx)\circ \pi_{x^{-1}gx}$ and\\
$\rho(\gx):A_{x^{-1}gx}\overset{\sim}{\to} A_g$ by isomorphisms.\\
In particular $\rho(\gx)(a_h)=\delta_{x^{-1}gx,h} \rho(\gx)(a_h)$
\end{enumerate}
\end{prop}
\begin{proof}

The equation (\ref{projectors}) means that the
$\rho(\elt{g}{e})$ act as projectors  and since that $\rho
(\elt{\mathbf 1}{e}) = id_A$ the first claim follows from  equation
(\ref{projectors}) by setting
\begin{equation}
A_g:=\rho(\elt{g}{e})(A)\end{equation}

For the first part of the second claim, we notice that
\begin{equation}\gx=
\elt{g}{e}\gx \;  \elt{(x^{-1}gx)}{e} \end{equation} which implies
the statement in conjunction with (1). For the second part, we
calculate that
$\elt{x^{-1}gx}{x^{-1}}\gx=\theta_{x^{-1}gx}(x^{-1},x)\elt{x^{-1}gx}{e}$
and since $\theta_{xgx^{-1}}(x^{-1},x)\neq 0$ and
$\rho(\elt{x^{-1}gx}{e})|_{A_{x^{-1}gx}}=\pi_{x^{-1}gx}|_{A_{x^{-1}gx}}=id$,
the claim follows.\end{proof}

\begin{nota}
\label{phinota} It will be convenient to denote $\rho(\elt{\mathbf
1}{g})$ by $\phi(g)$. For any $\dkg$ module $A$ we let
$A_g:=Im(\rho(\elt{g}{e}))$ and denote
the projection by $\pi_g$. Notice that then
\begin{equation}
\rho(\elt{g}{x})=\phi(x)\circ \pi_{x^{-1}gx}=
\phi(x)|_{A_{x^{-1}gx}}:A_{x^{-1}gx}\to A_g
\end{equation}
\end{nota}

\begin{rmk}
\label{dkgmodrem}

If $\beta\equiv 1$ then $\phi$ yields a $k[G]$ module structure on
$A$ while the grading corresponds to the $k[G]$ co-module structure
given by $a_g\mapsto a_g\otimes g$, moreover one can check that
these two structures are compatible so as to form a crossed $\dkg$
module in the sense of \cite{Kassel}, as is well known.

\end{rmk}

\subsection{DPR Induction}
\label{DPRindpar}
A very useful tool in the theory of the twisted Drinfel'd double is
the Dijkgraaf--Pasquier--Roche (DPR) induction \cite{DPR}.

For any $\a\in Z^2(G,\uk)$ let
 $R^{\a}(G)$ be the group of $\a$ twisted representations, that is maps
$\rho:G\to GL(V)$ with $\rho(g)\rho(h)=\a(g,h)\rho(g,h)$. We
write  $C(G)$ for the set of conjugacy classes of $G$.
With this notation DPR induction allows
one to constructively prove the following result.

\begin{thm}
\cite{DPR} \label{DPRthm} $Rep(\dbkg)\morita \bigoplus_{[g]\in C(G)}
R^{\theta_g}(Z(g))$.
\end{thm}

\begin{rmk}
One can also view the theorem above as following from the Morita equivalence
of the loop version of the inertia groupoid and the fiber product/stack
version of the inertia groupoid.
\end{rmk}

A very nice compilation of the results is given in \cite{ACM}. We also
 review the DPR induction process below.

\begin{rmk} We wish to point our several facts:
\begin{enumerate}
\item  Notice that the individual $R^{\a}(G)$ do not form rings. The product
induced by the tensor product on the
underlying modules is rather
from $R^{\a}(G)\otimes R^{\beta}(G)\to R^{\a\beta}(G)$. The direct
sum over the $\theta_g$ is in a certain
sense  ``closed'' under this
operation,  whence the product structure.
We refer to \cite{DPR} for the details, but also see \S \ref{DPRindpar} below.
\item The product in $Rep(\dbkg)$ is
not associative for general $\beta$, but only braided associative, with
the braiding given by the Drinfel'd associator $\Phi$.
See paragraph \S\ref{braidedpar}.

\item We write $k^{\a}[G]$ for the twisted group ring that
is $\bigoplus_{g\in G} k1_g$ with multiplication $1_g1_h=\a(g,h)1_{gh}$.
It is worth remarking that $G$ acts by conjugation
$\rho(g)1_h=\eps(g,h)1_{ghg^{-1}}$ with $\eps(g,h)=\frac{\a(g,h)}
{\a(ghg^{-1},g)}$. With this action (see e.g.\ \cite{Kar}):
\begin{equation}
(k^{\a}[G])^G\otimes {\mathbb C}=R^{\a}(G)\otimes {\mathbb C}
\end{equation}
Also, a module over $k^{\a}[G]$ is the same as an $\a$ twisted representation.
\end{enumerate}
Here and everywhere the superscipt $G$ denotes the $G$--invariants.
\end{rmk}

\begin{df}\cite{DPR}
\label{dprind}
Fix $\beta$ and $g\in G$.
Given $(V,\lambda)$ a left
 $\theta_g$ twisted representation of $Z(g)$ the {\it DPR induced}
representation is
$\dprind(V):= \field[G]\otimes_{\field^{\theta_g}[Z(g)]}V$
where for the tensor product $\field^{\theta_g}$ acts on the right on $\kg$
via
$x\rho(h)=\theta_{xgx^{-1}}(x,h)xh$
with the action of $\dbkg$ given by
\begin{equation}
\elt{h}{x}(r\otimes v):=\delta_{h,xrg(xr)^{-1}} \theta_h(x,r) xr\otimes v
\end{equation}
\end{df}

\begin{rmk}
Notice that if one chooses representatives $x_i$ for $G/Z(g)$ then the action
amounts to
\begin{eqnarray}
\elt{h}{x}(x_i\otimes v)&=& \delta_{h,xrg(xr)^{-1}}\theta_h(x,x_i)
xx_i \otimes v\nn\\
&=&
\delta_{h,xrg(xr)^{-1}}\theta_h(x,x_i) x_kz \otimes v\nn\\
&=&\delta_{h,xrg(xr)^{-1}}
\frac{\theta_h(x,x_i)}{\theta_{x_kgx_k^{-1}}(x_k,z)}
x_k\rho(z) \otimes v\nn\\
&=&\delta_{h,x_kgx_k^{-1}}
\frac{\theta_h(x,x_i)}{\theta_{h}(x_k,z)} x_k \otimes \lambda(z)(v)
\end{eqnarray}
which is the formula one can find for instance in \cite{ACM}.
\end{rmk}

\subsection{An exterior tensor product}

Recall \cite{wisc,orb} that for $G$--graded spaces
$A=\bigoplus_{g\in G}$ and $B=\bigoplus_{g\in G} B_g$ there is
another natural tensor product, which is given by
\begin{equation}
A\widehat\otimes B:= \bigoplus_{g\in G} A_g\otimes B_g
\end{equation}

\begin{prop}
If $A$ is a $\dkgv{\beta}$ module and $B$ is a $\dkgv{\beta'}$
module then $A\widehat\otimes B$ is a $\dkgv{\beta\beta'}$ module
via the diagonal action $\widehat\Delta(\gx)=\gx\otimes \gx$.
\end{prop}

\begin{proof}
First notice that indeed $A$ and $B$ are $G$--graded by Proposition
\ref{gradingprop}.

We need to check that
\begin{equation}
\widehat \Delta(\gx \hy)(a_k\otimes b_k)=\widehat \Delta
(\gx)(\widehat \Delta(\hy)(a_k\otimes b_k))
\end{equation}

For this to be non-zero, we need $h=yky^{-1}$ and $g=xyk(xy)^{-1}$,
so fix these values, then
$\gx\hy=\theta^{\omega}_{g}(x,y)\elt{g}{xy}$ in any $\dkgv{\omega}$.
Also, notice that by a simple substitution into the definitions
$\theta_g^{\beta\beta'}(x,y)=\theta^{\beta}_g(x,y)\theta^{\beta'}_g(x,y)$,
with which the claim follows. Finally
$\widehat\Delta(\elt{g}{e})=\elt{g}{e}\otimes \elt{g}{e}$ which
means that indeed the degree $g$ part of $A\widehat\otimes B$ is
given by $A_g\otimes B_g$.
\end{proof}

\subsection{A second exterior tensor product}
Notice that $\dkg$ as a vector space is actually bi--graded by
$G\times G$ and for bi--graded modules, there is again a tensor
product:

Now given any bi--graded $A=\bigoplus_{(g,x)\in G\times G} A_{g,x}$
and $B=\bigoplus_{(g,x)\in G\times G} B_{g,x}$ we define

\begin{equation}
A\widehat{\widehat\otimes} B:= \bigoplus_{(g,x) \in G\times G}
A_{g,x}\otimes B_{g,x}
\end{equation}

Of course this is just  $\hotimes$ for the group $H=G\times G$, but
since we consider the group $G$ to be fixed this notation will be very useful.

\begin{lem}
\label{hhotimeslem}
When using the diagonal product:
$\dbkg \hhotimes \dkgv{\b'}=\dkgv{\b\b'}$.
\end{lem}

\begin{proof}
Straightforward calculation.
\end{proof}

\section{$G$-Frobenius algebras}
\subsection{Frobenius algebras}
We wish to recall that there are two notions of Frobenius algebra.
The first goes back to Frobenius and is given as follows:

\begin{df}
A {\em Frobenius algebra} is a finite dimensional commutative
associative unital algebra $A$ together with a non--degenerate
symmetric pairing $\eta$ that is invariant, that is
\begin{equation}
\label{inveq} \eta(a,bc)=\eta(ab,c)
\end{equation}

A {\em possibly degenerate Frobenius algebra} is the same data as
above only that we do note require that $\eta$ is non--degenerate.
\end{df}

In the categorical setting there is the notion of a Frobenius
algebra object in a monoidal category.

\begin{df}
A {\it non--unital Frobenius algebra object} or Frobenius object for
short in a monoidal category ${\mathscr C}$ is an associative
commutative algebra object, which is also a co--associative
co--commutative object given by a datum $(A,\mu:A\otimes A\to A,
\Delta:A\to A\otimes A)$ that additionally satisfies

\begin{equation}
\label{frobeq} \Delta\circ \mu=(\mu\otimes id)\circ(id\otimes
\Delta)= (id\otimes \mu)\circ(\Delta\otimes id)
\end{equation}

A {\it Frobenius algebra object} is the data above
together with a unit $\unit: \unitc \to A$
 and a co-unit $\eps:A\to \unitc$, where
$\unitc$ is the unit object of $\C$.
\end{df}

\begin{rmk}
Notice that a Frobenius algebra always gives a Frobenius algebra
object in the monoidal category  ($k-Vect,\otimes$), by letting
$\Delta$ be the adjoint of $\mu$ with respect to the pairing. The
co--unit is given by pairing with the unit of the algebra.

Vice--versa if $A$ is a Frobenius algebra object in
($k-Vect,\otimes$) then $A$ with its unit, multiplication and
$\eta(a,b):=\eps(\mu(a\otimes b))$ is a possibly degenerate
Frobenius algebra.
\end{rmk}

\subsection{F-Traces and Trace Elements}
One main difference between the finite dimensional and the
non--finite dimensional case is the existence of traces.
In the finite dimensional case, for any operator $\phi\in Aut(A)$ we
can consider $Tr(\phi)$.
The trace actually has an analog in the Frobenius object case, for
this we need an expression in terms of the morphisms.

\begin{prop} For a Frobenius algebra, let $1_k$ be the unit in $k$ then
\begin{equation}
Tr(\phi)=\eps(\mu(\phi\otimes id) \Delta(\unit(1_k)))
\end{equation}
\end{prop}

\begin{proof}
Let $1_A=\unit(1_k)$ be the unit of $A$ and let $\Delta_i$ be a
basis of $A$. If $g_{ij}=\eta(\Delta_i,\Delta_j)$ is the metric and
$g^{ij}$ is its inverse then
$\Delta(\unit(1_k))=\Delta(1_A)=\sum_{ij} g^{ij}\Delta_i\otimes
\Delta_j$, since
$$
\eta\otimes \eta (\Delta_k\otimes \Delta_l,\Delta(1_A)):=
\eta(\Delta_k\Delta_l,1_a)=\eta(\Delta_k,\Delta_l)=g_{kl}
$$
and
$$
\eta\otimes \eta (\Delta_k\otimes \Delta_l,\sum_{ij}
g^{ij}\Delta_i\otimes \Delta_j)= \sum_{ij} g_{ki}g^{ij}g_{jl}
=\sum_i g_{ki}\delta_{i,l}=\delta_{k,l}
$$
and hence if $\check\Delta_i:=\sum_jg^{ij}\Delta_j$ is the inverse
basis
\begin{multline*}
\eps(\mu(\phi\otimes id) \Delta(\unit(1_k)))=\eps(\sum_{ij}g^{ij}
\phi(\Delta_i) \Delta_j)=\\\sum_{i,j}\eta(\sum_{ij}g^{ij}
\phi(\Delta_i),
\Delta_j)=\sum_i\check\Delta_i(\phi(\Delta_i))=Tr(\phi)
\end{multline*}
\end{proof}

\begin{df}
Given a Frobenius algebra object $A$ and $\phi \in Aut(A)$ we define
the {\it F-Trace} $\tau(\phi): \unitc \to \unitc$ of $\phi$ via
\begin{equation}
\tau(\phi):=\eps\circ \mu \circ (\phi\otimes id) \circ \Delta \circ \eta
\end{equation}
\end{df}

\begin{rmk}
If $\unitc=k$ and all morphisms
are $k$--linear the map $\tau(\phi)$ is of course given by its
value on $1_k$. In this case we will not distinguish between the map and
this value.
\end{rmk}

\begin{prop}
\label{monoidalprop} Let $\mathcal{F}$ be a monoidal functor with
values in vector spaces for a category with products given by the
monoidal structure. Also assume that $\mathcal F$ has pull--backs,
push--forwards and satisfies the projection formula for the diagonal
morphisms. Then for any object $V$ it gives rise to a Frobenius
algebra object and hence F-traces.

\end{prop}
\begin{proof}
We let $\mu$ be given by the pull--back along the diagonal $\Delta_V:
V\to V\times V$ where the co--multiplication is given by
push--forward along the diagonal:
$\mu=\Delta_V^*,\Delta=\Delta_{V*}$.

The equation (\ref{frobeq}) is guaranteed by the projection formula.
On one hand:
\begin{equation}
\Delta_{V*}(\Delta_V^*(\F_1 \otimes \F_2))=(\F_1\otimes
\F_2)\Delta_{V*}(1)
\end{equation}
On the other hand:
\begin{multline}
 (\Delta_{V}^*\otimes id )(id\otimes \Delta_{V*})(\F_1\otimes
\F_2) =(\Delta_{V}^*\otimes id )(\F_1 \otimes
\Delta_{V*}(\Delta_V^*(1 \otimes \F_2)))\\
= (\Delta_{V}^*\otimes id) (\F_1 \otimes (1 \otimes
\F_2)\Delta_{V*}(1))=\sum \F_1 \Delta^{(1)}\otimes \F_2
\Delta^{(2)}=(\F_1\otimes \F_2)\Delta_{V*}(1)
\end{multline}
where we used Sweedler's notation $\Delta_{V*}(1)=\sum
\Delta^{(1)}\otimes \Delta^{(2)}$ and analogously for the third
equation.

The co--unit is  furnished by the push--forward to the unit of the
monoidal category  which is a final object, and the unit of the
Frobenius algebra object by the pull--back from it. In our cases of
interest this will be a point or the one--dimensional vector space
of the ground field.

\end{proof}

\begin{cor}
In the situation above, we also obtain pre-Frobenius algebras in the
sense of \cite{JKK2}, where the trace element is the morphism given
by $\forall a\in A: a\mapsto \tau(\lambda_a)$ that is the F-trace of
the morphism of left multiplication by $a$; $\lambda_a(b):=ab$.
\end{cor}

\begin{ex}
Notice that this gives the canonical trace elements considered in
\cite{JKK2} for the pre--Frobenius algebras $A^*(V)$ and $K_0(V)$,
which are prime examples of Frobenius algebra objects, that give
rise to possibly degenerate Frobenius algebras, as they might be
infinite dimensional. Here $\eps = \int$ or $\chi$ respectively,
which are the push--forwards to a point.
 For example in $A^*$,
 we can calculate the F--trace $\l_v$ ---which is the operation of
left multiplication by  $v$--- to be given by
\begin{eqnarray}
\tau(\lambda_v)&=&\int_V \Delta_V^*[(v\otimes 1_V)\cup(\Delta_{V_*}(1))]\nn\\
&=&\int_V (v \cup \Delta^{(1)})\cup \Delta^{(2)}\nn\\
&=&
\int_V v \cup \Delta_V^*\Delta_{V_*}(1)\nn\\&=&\int_V v \cup c_{\rm top}(TV)\\
\end{eqnarray}
where we used the notation of the last Proposition for the co--product.
This is exactly the expression appearing in \cite{JKK2}. The
analogous statement of course holds for K--theory.
\end{ex}

\subsection{Twisted Frobenius objects}
In general there is a twisted version of Frobenius algebra objects.
This appears in the definition of $G$--Frobenius algebras and is necessary
for considerations concerning singularities with symmetries, see e.g.\ \cite{wisc,orb,singorb}.  We again fix a monoidal category $\C$.

\begin{df}
Let $\twist$ be an even invertible element in $\C$.

A $\twist$--twisted
Frobenius algebra object is the datum
$(A,\mu: A\to A\otimes A, \Delta:A\to A\otimes A\otimes
\twist^{\otimes 2},\unit:\unitc \to A, \counit:A\to \twist^{\otimes 2})$
such that (\ref{frobeq}) is satisfied,
where $\mu$ is associative commutative, $\counit$ is co--associative,
co--commutative,
 $\unit$ is a unit, and $\counit$ is a
co--unit using the isomorphism $m:\twist \otimes \twist^{-1}\iso
\unitc$. More precisely:

$$\begin{CD}
 \twist^{\otimes -2} \otimes A \otimes \twist^{\otimes 2}
&@<\counit \otimes Id\otimes Id<<&
A \otimes A \otimes \twist^{\otimes 2} &@>id \otimes
\counit>>&A\otimes \twist^{\otimes -2}  \otimes  \twist^{\otimes 2}\\
@Vm\circ m^{\otimes 2}VV&&@A\Delta AA&&@VVm\circ m^{\otimes 2}V\\
\unitc\otimes A&@>>>&A&@<<<&A \otimes \unitc
\end{CD}
$$
where on the left $m^{\otimes 2}$ is $m$ applied to the 1st and 4th and the
2nd and 5th component and then to the two copies of $\unitc$
and  on the right to the  2nd and 4th and to the 3rd and 5th and then
again to the two copies of $\unitc$.
\end{df}

\begin{rmk}
One could of course twist $A\leadsto \bar A:=A\otimes \twist^{-1}$ and
obtain similar operations and axioms. In the language of \cite{wisc,orb}
this is the Ramond twist or Ramond sector.
\end{rmk}

\subsection{$G$-Frobenius algebras}
First we recall the main definition see \cite{wisc,orb}:
\begin{df}
A {\it $G$--Frobenius algebra} or GFA for short,
over a field  $k$ of characteristic 0 is
$<G,A,\circ,1,\eta,\varphi,\chi>$, where

\begin{tabular}{ll}
$G$&finite group\\
$A$&finite dim $G$-graded $k$--vector space \\
&$A=\oplus_{g \in G}A_{g}$\\
&$A_{e}$ is called the untwisted sector and \\
&the $A_{g}$ for $g \neq
e$ are called the twisted sectors.\\
$\circ$&a multiplication on $A$ which respects the grading:\\
&$\circ:A_g \otimes A_h \rightarrow A_{gh}$\\
$1$&a fixed element in $A_{e}$--the unit\\
$\eta$&non-degenerate bilinear form\\
&which respects grading i.e. $g|_{A_{g}\otimes A_{h}}=0$ unless
$gh=e$.\\
\end{tabular}

\begin{tabular}{ll}
$\varphi$&an action  of $G$ on $A$
(which will be  by algebra automorphisms), \\
&$\varphi\in \mathrm{Hom}(G,\mathrm{Aut}(A))$, s.t.\
$\varphi_{g}(A_{h})\subset A_{ghg^{-1}}$\\
$\chi$&a character $\chi \in \mathrm {Hom}(G,k^{*})$ \\

\end{tabular}

\vskip 0.3cm

\noindent Satisfying the following axioms:

\noindent{\sc Notation:} We use a subscript on an element of $A$ to
signify that it has homogeneous group degree  --e.g.\ $a_g$ means
$a_g \in A_g$-- and we write $\varphi_{g}:= \varphi(g)$ and
$\chi_{g}:= \chi(g)$.

\begin{itemize}

\item[a)] {\em Associativity}

$(a_{g}\circ a_{h}) \circ a_{k} =a_{g}\circ (a_{h} \circ a_{k})$
\item[b)] {\em Twisted commutativity}

$a_{g}\circ a_{h} = \varphi_{g}(a_{h})\circ a_{g}$
\item[c)]
{\em $G$ Invariant Unit}:

$1 \circ a_{g} = a_{g}\circ 1 = a_g$

and

$\varphi_g(1)=1$
\item[d)]
{\em Invariance of the metric}:

$\eta(a_{g},a_{h}\circ a_{k}) = \eta(a_{g}\circ a_{h},a_{k})$

\item[i)]
{\em Projective self--invariance of the twisted sectors}

$\varphi_{g}|A_{g}=\chi_{g}^{-1}id$

\item[ii)]
{\em $G$--Invariance of the multiplication}

$\varphi_{k}(a_{g}\circ a_{h}) = \varphi_{k}(a_{g})\circ
\varphi_{k}(a_{h})$

\item[iii)]

{\em Projective $G$--invariance of the metric}

$\varphi_{g}^{*}(\eta) = \chi_{g}^{2}\eta$

\item[iv)]
{\em Projective trace axiom}

$\forall c \in A_{[g,h]}$ and $l_c$ left multiplication by $c$:

$\chi_{h}\mathrm {Tr} (l_c  \varphi_{h}|_{A_{g}})=
\chi_{g^{-1}}\mathrm  {Tr}(  \varphi_{g^{-1}} l_c|_{A_{h}})$
\end{itemize}

We call a $G$--Frobenius algebra {\em strict}, if $\chi\equiv 1$.
\end{df}

\begin{rmk}
It was shown in \cite{disc} that a GFA is a module over $\dkg$
and moreover proved that it is a $k[G]$ module algebra and a $k[G]$ co--module
algebra. The first part also follows from Remark \ref{dkgmodrem}.
\end{rmk}

\begin{ex}
Important examples are furnished by the twisted group rings
$k^{\a}[G]$ with $\a\in Z^2(G,\uk)$. This group actually acts
on the set (category) of $G$--Frobenius algebras
 through $\hotimes$ and gives rise to the action of
discrete torsion, see  \cite{disc} for full details.
\end{ex}

\begin{prop}
A $G$--Frobenius algebra with character $\chi$
 is a unital, associative, commutative
algebra object in the category $\dkg$--mod. It moreover defines a $k_{\chi}$
twisted Frobenius algebra object, where $k_{\chi}$ is the 1--dimensional
$\dkg$--module concentrated in group degree $e$ with $G$ action on $k$
given by the character $\chi$.
\end{prop}

\begin{proof}
This follows in a straightforward fashion, by reinterpreting the
pertinent diagrams using the braided monoidal structure.

Since $\beta\equiv 1$ associativity in the category $\dkg$--{\it
Mod} is just the ordinary associativity a).

Let $\mu$ denote the multiplication in $A$. In view of equation
(\ref{action}) the $G$--invariance of the multiplication ii) is
equivalent to $\mu:A\otimes A \to A$ being a morphism in the
category $\dkg$--{\it Mod}.

Using the equation (\ref{braid}) we see that the condition that the
following diagram commutes ---which is the commutativity in
$\dkg$--{\it Mod}--- is equivalent to the condition b) of twisted
commutativity.

$$\begin{CD}
A\otimes A &@>\mu>>& A\\
@V c_{A,A} VV&&@VV id V\\
A\otimes A &@>>>& A\\
\end{CD}
$$

The fact that the unit is invariant is equivalent to the diagram
$$\begin{CD}
k\otimes A &@>\unit \otimes Id>>A \otimes A &@<id \otimes
\unit<<&A\otimes k\\
&\nwarrow&&@V\mu VV&\nearrow&\\
&&&A&&&
\end{CD}
$$
being a diagram of $\dkg$ modules where $k$ has the structure of a
trivial $\dkg$ module.

 We define the co--unit via
$\counit(a):=\eta(a,1_k)$. Then the projective $G$--invarince
of the metric iii) becomes the condition on the co--unit in a twisted
Frobenius algebra.

We set $\Delta:=\mu^{\dagger}$, that is the adjoint of the
multiplication under the non--degenrate metric $\eta$. Then the
invariace of the metric d) together with the projective
$G$--invariance iii) yields the Frobenius equation (\ref{frobeq}).

\end{proof}

\begin{thm}\label{frobthm}
A $G$--Frobenius algebra with character $\chi$ is precisely a
$\twist$--twisted Frobenius algebra object $\dkg$--{\it Mod} with
the following additional restrictions
\begin{itemize}
\item[1)] The associated pairing $\eta=\counit\circ \mu$
is non--degenerate.
\item[2)] Denoting the $\dkg$ action induced
by the $G$ action $\varphi$  by $\rho$ the following two axioms hold
\begin{itemize}
\item[(T)]
$\rho(v^{-1})=\chi^{-1}$ for a character $\chi\in \Hom(G,\uk)$
\item[(S)] Using the Notation \ref{phinota} let
 $l_c$ denote the left multiplication by $c$:
$c\in A$:
\begin{equation}
\chi_h\tau (l_c\circ \rho(\elt{hgh^{-1}}{h}))=
\chi_{g^{-1}}\tau(\rho(\elt{h}{g^{-1}})\circ l_c)
\end{equation}
where $\tau$ is the $F$--trace.
\end{itemize}

\end{itemize}

\end{thm}
\begin{proof}
Given a GFA, it is a  unital, associative, commutative algebra
object in $\dkg$--{\it Mod} by the above Proposition and it also
satisfies the additional axioms. By the proposition
\ref{gradingprop}, we see that any $\dkg$--{\it Mod} is $G$ graded
and has an action of $G$ by automorphisms of $G$ given by $\phi$ of
Notation \ref{phinota}, which act in the prescribed way. Now by the
proof of the Proposition above, we have that a unital associative
commutative algebra object satisfies the axioms a),b),c),ii). What
remains to be shown is that the multiplication preserves the
$G$--grading, but this follows from the fact that
$\Delta(\elt{k}{e})=\sum_{gh=k}\elt{g}{e}\otimes \elt{h}{e}$ so that
if the multiplication is a morphism  the multiplication is graded
since the $\elt{g}{e}$ act as projectors. Explicitly
$$\rho(\elt{k}{e})(a_gb_h)=\mu\circ (\rho\otimes\rho)(
\Delta(\elt{k}{e}))(a_g\otimes a_h)=\delta_{k,gh}a_gb_h$$

It is clear that $\eta=\counit\circ \mu$ defines
a pairing given a Frobenius algebra object and
as above vice--versa defines $\counit$ in the presence of a unit.
 The invariance of the metric d) follows from the Frobenius
equation and the structure of the co--unit. The latter is also equivalent
to the projective $G$--invariance of the metric iii).

For the equivalence of the
projective trace axiom with (S), we recall that the elements $\elt{g}{x}$ act
as explained in Notation \ref{phinota}. Notice the if $c\notin A_{[gh]}$ then
both sides are zero.
In the same notation with the definition of $v^{-1}$ given equation (\ref{veq})
condition (T) is just condition i).

\end{proof}

Here S and T stand for the generators of $SL(2,{\mathbb Z})$ and are
a reminder that these axioms correspond to the invariance of the
conformal blocks.\footnote{As someone suggested,
 $S$ could of course also stand for {\it Spur}.}

Dropping the condition 1) we come to the main definition of the paragraph.

\begin{df}
We define a {\it  $G$--Frobenius algebra object} to be a Frobenius algebra
object in the category $\dkg-Mod$ which satisfies the axioms (S) and (T).
\end{df}

\begin{rmk}
Going beyond the aesthetics and the practicality of the above
definition, it is a necessary generalization if we are to deal with
the stringy Chow ring or Grothendieck $K$--theory of a global
quotient stack as in \cite{JKK2}, where the natural metric may be
degenerate. See the next paragraph for details.
\end{rmk}

\subsection{The Drinfel'd double as a $G$--Frobenius algebra}
\label{dprindsection}
We have seen that any GFA is actually a $\dkg$ module. Now as it
happens $\dkg$ is itself a $\dkg$ module, but not quite a $G$--Frobenius
algebra for general $G$.
This is because the $G$--degree of $\elt{g}{x}$ is $g$
and the multiplication is not multiplicative in $g$ but rather in $x$.

Notice that the elements $\elt{g}{x}$ with $[g,x]=e$ form a
subalgebra $\subdbkg$ of $\dbkg$ which is actually additively
isomorphic to $\bigoplus_{g\in G} \field^{\theta_g}(Z(g))$. In the
case that $G$ is Abelian of course $\subdbkg=\dbkg$.

\begin{prop}
$\subdbkg$ is a GFA for the $\dkg$ action given by
\begin{equation}
\rho(\elt{g}{x})(\elt{h}{y})=
\frac{\theta_{xhx^{-1}}(x,y)}{\theta_{xhx^{-1}}(xyx^{-1},x)}\delta_{g,xyx^{-1}}\elt{xhx^{-1}}{xyx^{-1}}
\end{equation}
which means that the $G$--degree of $\elt{h}{y}$ is $y$.
\end{prop}

\begin{proof}
This follows from the fact that the each $k^{\theta_g}[Z(g)]$ is
actually a $Z(g)$--FA. This means for instance that it satisfies all
the axioms for the $Z(g)$ action pertaining to the $Z(g)$ alone. The
other axioms then follow from the $G$--equivariance of the
$\theta_g$ or are straightforward. For $\beta=1$ the statement also
follows from Proposition \ref{korbinertia}.
\end{proof}

\begin{rmk}
In the case of $\dkg$ if one uses the grading that the $G$ degree of
$\elt{g}{x}$ is $x$ so that the multiplication is indeed
$G$--graded, then  twisted commutativity dictates that
$\rho(\elt{{\bf
1}}{h})(\elt{g}{x})=\elt{hxgx^{-1}h^{-1}}{hxh^{-1}}$. In turn
postulating the compatibility of this $G$ action with the
multiplication requires that $[g,x]=e$.
\end{rmk}

\begin{df}
We call a GFA a {\it free GFA}, if it is of the form
$A=\field^{\theta_g}[G]\otimes A_e$ for a Frobenius algebra $A_e$
that is a $G$--module, with the multiplication given by the diagonal
multiplication, the $G$--degree of $g\otimes a$ being $g$ and the
$G$--action given by the conjugation action on the left factor and
the postulated $G$ action on $A_e$.
\end{df}

\begin{rmk}
Notice that in this case we have a second $\field^{\theta_g}[G]$
action, given by multiplication from the left on the factor
$\field^{\theta_g}[G]$. This action sends $\lambda_h:A_g\to A_{hg}$.
This is similar to the quantum symmetry considered in \cite{disc}.
\end{rmk}

\begin{df}
\label{dpralgebrainddef} Given $\b\in Z^3(G,\uk)$, if
$A=\bigoplus_{g \in G} B_g$ is the direct sum of free
$Z(g)$--Frobenius algebras $B_g=\field^{\theta_g}[Z(g)]\otimes B_e$
then we define the {\it DPR induced free} algebra
$\dprind(A):=\bigoplus \field[G]\otimes_{\field^{\theta_g}[Z(g)]}
B_g\iso \dbkg\otimes B_e$, where the action is analogous to
Definition \ref{dprind} and the algebra structure is the diagonal
algebra structure.
\end{df}

\begin{rmk}
At the moment we do not see how to induce this algebra in the
non--free case. Geometrically this amounts to the fact that on the
inertia, the automorphisms have to commute so that the double
twisted sectors for non--commuting elements are not accessible. Also
in the general case, the double twisted sectors $A_{x^{-1}gx,x}$ for
$x\in G$ are not equidimensional. This is however an interesting
detail which should be studied further, but is unfortunately beyond
the scope of the present considerations.
\end{rmk}

\section{Orbifold cohomology and K-theory}
In this section, we recall the various stringy functors introduced
in \cite{JKK2} and re-express them in the current framework. First,
we recall from \cite{JKK2} that we have the following stringy
functors for a global quotient $(X,G)$, $\calF\in\{A^*, H^*,K_0,
K^{top}\}$ as well as isomorphisms
 $\Ch: K_0(X,G) \to A^*(K,G)$ and $\Ch: K^{top}(X,G)\to H^*(X,G)$.
Then we also recall the stack versions of these functors and maps for
a suitably nice stack $\cx$. In order to simplify things we will
work over ${\mathbb Q}$ or extensions of it, see however Remark \ref{zrem}.

\subsection{General setup -- global quotient case}
We recall the setup as in the global part of \cite{JKK2}. We
simultaneously treat two flavors of geometry: algebraic and
differential. For the latter, we consider a stably almost complex
manifold $X$ with the action of a finite group $G$ such that the
stably almost complex bundle is $G$ equivariant. While for the
former $X$ is taken to be a smooth projective variety with a $G$--action.

In both situations for $m\in G$ we denote the fixed point set of $m$
by $X^m$ and let
\begin{equation}
I(X)=\amalg_{m\in G} X^m
\end{equation}
be the inertia variety.

We let ${\mathcal F}$ be any of the functors $H^*,K_0,A^*,K^{\rm
top}$, that is cohomology, Grothendieck $K_0$, Chow ring or
topological $K$--theory with ${\mathbb Q}$ coefficients,  and define
\begin{equation}
\F_{stringy}(X,G):= \F(I(X))= \bigoplus_{m\in G} \F(X^m)
\end{equation}
{\em additively}.

We furthermore set
\begin{equation}
\eu{\F}(E)=\begin{cases} c_{top}(E)&\text{ if } \F=H^* \text { or }
A^*
\text{ and $E$ is a bundle}\\
 \lambda_{-1}(E^*)&\text{ if } \F=K \text { or } K^{top}\\
 \end{cases}
\end{equation}

Notice that on bundles $\e$ is multiplicative. For general
$K$--theory elements we set
\begin{equation}
\eut{\F}(E)=\begin{cases} c_{t}(E)&\text{ if } \F=H^* \text { or } A^*\\
 \lambda_{t}(E^*)&\text{ if } \F=K \text { or } K^{top}
\end{cases}
\end{equation}

\subsection{The stringy product}
For $m \in G$ we let $\Xm$ be the fixed point set of $m$ and for a
triple $\bm=(m_1,m_2,m_3)$ (or more generally an n--tuple)
such that  $\prod m_i={\bf 1}$ (where
${\bf 1}$ is the identity of $G$) we let $X^{\bm}$ be the common
fixed point set, that is the set fixed under the subgroup generated
by them.

In this situation, recall the following definitions. Fix $m \in G$
let $r=ord(m)$ be its order. Furthermore let $W_{m,k}$ be the
sub--bundle of $TX|_{X^{\bm}}$ on which $m$ acts with character
$\exp(2\pi i \frac{k}{r})$, then
\begin{equation}
S_m=\bigoplus_k \frac{k}{r} W_{m,k}
\end{equation}
Notice this formula is invariant under stabilization.

We also wish to point out that using the identification
$X^m=X^{m^{-1}}$
\begin{equation}
\label{normaleq} S_m\oplus(S_{m^{-1}})=N_{X^m/X}
\end{equation}
where for an embedding $X\to Y$ we will use the notation $N_{X/Y}$
for the normal bundle.

Recall from \cite{JKK2} that in such a situation there is a product
on $\F(X,G)$ which is given by

\begin{equation}
\label{proddefeq}
 v_{m_1}*v_{m_2}:=\check e_{m_3*}(e^*_1(v_{m_1})e_2^*(v_{m_2})\e(\R(\bm)))
\end{equation}

where the obstruction bundle ${\mathcal R}(\bm)$ is defined by
\begin{equation}
{\mathcal R}(\bm)=  S_{m_1}\oplus S_{m_2}\oplus S_{m_3} \ominus
N_{X^{\bm}/X}
\end{equation}
and the $e_i:X^{m_i}\to X$ and $\check e_3:X^{m_3^{-1}}\to X$ are
the inclusions. Notice, that as it is written $\R(\bm)$ only has to
be an element of K-theory with rational coefficients, but is
actually indeed represented by a bundle \cite{JKK2}.

\begin{rmk}
\label{zrem}
This bundle and hence
the multiplication below are actually defined over ${\mathbb Z}$.
The point is that in \cite{JKK2} we identified $\R(\bm)$ as a bundle
and true representation in the representation ring. Since there
is no torsion in this ring the bundle is identified over   ${\mathbb Z}$.
\end{rmk}

\begin{rmk}
The first appearance of a push--pull formula was given in \cite{CR}
in terms of a moduli space of maps. The product was for the $G$
invariants, that is for the $H^*$ of the inertia orbifold and is
known as Chen--Ruan cohomology. In \cite{FG} the obstruction bundle
was given using Galois covers establishing a product for $H^*$ on
the inertia variety level, i.e.\ a $G$--Frobenius algebra as defined
in \cite{wisc,orb}, which is commonly referred to as the
Fantechi--G\"ottsche ring. In \cite{JKK1}, we put this global
structure back into a moduli space setting and proved the trace
axiom. The multiplication on the Chow ring $A^*$ for the inertia
stack was defined in \cite{AGV}. The representation of the
obstruction bundle in terms of the $S_m$ and hence the passing to
the differentiable setting as well as the two flavors of $K$--theory
stem from \cite{JKK2}.
\end{rmk}

The following is the key diagram:

\begin{equation}
\label{maindiagrameq}
\begin{matrix}
X^{m_1}&X^{m_2}&X^{m_3^{-1}}\\
e_1\nwarrow &\uparrow e_2 &\nearrow \check e_3\\
&X^{\bm}&
\end{matrix}
\end{equation}

Here we used the notation of \cite{JKK2}, where $e_3:X^{\bm}\to
X^{m_3}$ and $i_3: X^{m_3}\to X$ are the inclusion, $\vee: I(X)\to
I(X)$ is the involution which sends the component $X^m$ to
$X^{m^{-1}}$ using the identity map and $\check \imath_3= i_3\circ
\vee$, $\check e_3=\vee\circ e_3$. This is short hand notation for
the general notation of the inclusion maps $i_m:X^m\to X$, $\check
\imath_m:=i_m\circ \vee=i_{m^{-1}}$.

\begin{thm}\label{frobobthm}
The cases in which
$\F$ equals $H^*$ and $K^{top}$ yield $G$--Frobenius algebras.
In the cases of $A^*$ and $K_0$ the stringy functors are still
$G$--Frobenius algebra objects.
The co-multiplication is given by
\begin{equation}
\Delta(\F_{m_3})=
\sum_{m_1,m_2, m_1m_2=m_3}(\check e_{1*}\otimes \check e_{2*} )
\Delta_{X^{\bm}*}(e^*_3(\F_3)\e(\R(\bm))
\end{equation}
where $\Delta_{\Xm}: \Xm\to \Xm\times \Xm$ is that diagonal map.

The F-traces  $\tau(\lambda_c \phi_g,h):=\tau(\lambda_c\circ
\rho(\elt{ghg^{-1}}{g}))$ give the trace elements which were part of
the definition of pre--Frobenius algebra structures defined in
\cite{JKK2}.
\end{thm}

\begin{proof}
The first part about $H^*$ and $K^{top}$ is contained in
\cite{JKK2}. For $A^*$ or $K_0$  the verification of the Frobenius
condition (\ref{frobeq}) is somewhat tedious but straightforward
using analogous arguments as in Proposition \ref{monoidalprop}. We
will calculate the trace elements. We fix $a,b\in G$ and
$v_{[a,b]}\in A_{[a,b]}$ and calculate the F-trace
$\tau(\lambda_{v_{[a,b]}}\phi_b,a)$. For this we will need to set up
some notation and recall some results from \cite{JKK2}. We will use
the following notation analogous to \cite{JKK2}
$\bm'=([a,b],bab^{-1},a)$, $H':=<[a,b],bab^{-1}>\subset H:=<a,b>$.
We will also need the commutative diagram

$$\begin{CD}
X^H &@>j_2'>>& X^{H^{\prime}}\\
@V j_1' VV&&@VV\Delta_2' V\\
X^a&@>\Delta_1' >>& X^{bab^{-1}}\times X^a\\
\end{CD}
$$
where $j_1'$ and $j_2'$ are the inclusion morphisms, $\Delta_2'$ is the
diagonal map, and $\Delta_1'$ is the composition

$$
\begin{CD}
X^a @> \Delta_{X^a}>>
X^a\times X^a @>\phi(b)\times \vee >>
X^{bab^{-1}}\times X^{a^{-1}}
\end{CD}
$$
We denote the excess intersection bundle by $\mathscr{E}'$.
Also, we  recall that for a triple product
$v_{m_1}*v_{m_2}*v_{m_3}$ we have a special formula which actually is
the reason for associativity.

Let $\bm=<m_1,m_2,m_3, m_4=(m_1m_2m_3)^{-1}>$, and
$\bm'=<m_1,m_2,(m_1m_2)^{-1}>$
$X^{H^{\prime}}:=X^{\bm}$ and as usual  let
$e_i:X^{H^{\prime}} \to X^{m_i}$ be the inclusions and $\check{e_i}=\vee\circ e_i$,
then we have

$$
v_{m_1}*v_{m_2}*v_{m_3}=\check e_{m_4*}[(\prod e_i^*(v_{m_i}) \etop(\R(\bm))]
$$
where $\R(\bm)=\bigoplus S_{m_i}\ominus N(X^H)$.

Let $p_V:V\to pt$ be the projection to a point. In our case
$\bm=([a,b],bab^{-1},a^{-1},e)$ and $\bm'=([a,b],bab^{-1},a^{-1})$
and $H=\langle a, b\rangle$. Let $1_V$ be the unit in $\F(V)$. Then:

$$\Delta(1_X)=\sum_h e_{h*}\otimes e_{h^{-1}*}(e^*_h(1_X)
\etop( \R((h,h^{-1},e)))= \sum_h (id \otimes \vee)\Delta_{X^h}(1_{X^h})$$

So that the bi--degree $(h,h^{-1})$--part is just given by  $\ (id \otimes \vee)\Delta_{X^h}(1_{X^h})$

\begin{eqnarray}
\tau(\phi(b),a)&=&p_{X*}[\check e_{m_4*}[e_1^*(v_{[a,b]}\Delta_2'
(\Delta_1'(1))\etop(\R(\bm')) ]]\nn\\
&=&p_{X^H*}[e_1^*(v_{[a,b]}j'_{2*}(j_1^{\prime *}(1_{X^h})
\etop(\mathscr{E}')
\etop(\R(\bm')) ]\nn
\\
&=&p_{X^{H'}*}
[v_{[a,b]}|_{X^{H^{\prime}}} \etop(\mathscr{E}')j_2^{\prime *}
(\etop(\R(\bm'))) ]\nn\\
&=&p_{X^{H'}*}
[v_{[a,b]}|_{X^{H^{\prime}}} \etop(\mathscr{E}')\oplus j^{\prime *}_2
( \R(\bm'))) ]\nn\\
&=&p_{X^{H'}*}[v_{[a,b]}|_{X^{H^{\prime}}}
\etop(TX^H\oplus S_{[a,b]}|_{X^H}) ]\
\end{eqnarray}
which is the expression of \cite{JKK2}.
Here the last equality  follows from  the equality of the bundles
$ \mathscr{E}'\oplus j^{\prime *}_2( \R(\bm'))=
TX^H\oplus S_{[a,b]}|_{X^H} $ which fittingly was proved in
\cite{JKK2} (Theorem 5.5).

The traces $\tau(\lambda_v \phi_b,a)$ will of course
 be zero if $v$ is of pure $G$--degree different from $[a,b]$.

\end{proof}

\begin{prop}
\label{kunnethprop} Given $(X,G)$ and $(Y,G)$, $X\times Y$ has a
diagonal $G$ action and $\Forb((X\times Y,G))=\Forb((X,G))\widehat
\otimes \Forb((Y,G))$ where $\Forb$ is the global stringy version of
any of the functors $\F$ as defined in \cite{JKK2}.
\end{prop}

\begin{proof}
Straightforward by the K\"{u}nneth formula or relevant versions thereof.
\end{proof}

\subsection{The stack case}
In \cite{JKK2} a version of stringy K--theory or Chow
for general stacks was developed
as well. The important thing about the stringy K--theory in this case, which was also called full orbifold K--theory is that is it usually bigger than the global K-theory. In particular for a stack $\mathfrak{X}$ it was defined that
$K_{full}(\cx):=K(\ccx)$ where
$\ccx$ is the inertia stack. For a global quotient stack
we also defined
$K_{small}([X/G]):=K_{global}(X,G)^G$. Notice that this is
actually presentation
independent \cite{JKK2}.

In particular for a global quotient three theories where introduced
which are {\em additively} over given $\mathbb C$ as follows.

\begin{eqnarray}
K_{global}((X,G)):=&K(I(X,G)))&\iso\bigoplus_{g\in G}K(X^g)\\
K_{full}([X/G]):=&K(\mathfrak {I}[X/G])&\iso \bigoplus_{[g]} K([X^g/Z(g)])\\
K_{small}([X/G]):=&K_{global}((X,G))^G&\iso\bigoplus_{[g]}
K(X^g)^{Z(g)}
\end{eqnarray}
Here these are only {\em linear isomorphism} and the product is the one
given by the push--pull formula \ref{proddefeq}.
Notice they are all different.
It is however the case that $K_{small}$ is a subring of $K_{full}$
(see \cite{JKK2}).

\subsection{Comparing the different constructions
in the case of a global quotient} As mentioned above for global
quotient stacks we have $\Ksmall([X/G])\iso K(X,G)^G$ which is
isomorphic to $A^*$ or $H^*$, but also we have $\Kfull([X/G])$,
which is usually much bigger. Notice that  $\Ksmall(I(X,G),G)$ and
$\Kfull[(X/G)]$ are of the same size but have different
multiplications that is they {\em are additively isomorphic, but not
multiplicatively}.

\begin{prop}
\label{inertiacalcprop}
Additively:
\begin{eqnarray}
\label{korbeq}
\Korb(I(X,G),G)&=&K(\amalg_{x\in G} (\amalg_{g\in G} X^g)^x\\
&=& \bigoplus_{g\in G, x\in Z(g)} K(X^{\la
g,x\ra})\\
\end{eqnarray}
and for $\prod x_i=1$, and $g: x\in Z(g), h: y\in Z(h)$ the
multiplication is given by
\begin{multline}
\F_{g,x_1}*\F_{h,x_2}=\check
e_{x_3*}(e_{x_1}^*(\F_{g,x_1})e_{x_2}^*(\F{h,x_2})\R((x_1,x_2,(x_1x_2)^{-1}))\\
=\delta_{g,h}\; \F_{g,x_1}*_g\F_{g,x_2}
\end{multline}
where $*_g$ is the multiplication on $\Korb(X^g,Z(g))$, that is as
rings
\begin{equation}
\Korb(I(X,G),G)=\bigoplus_{g\in G} \Korb(X^g,Z(g))
\end{equation}
\end{prop}

\begin{proof}
Notice that if $g\neq h$ then the pull--backs land in different
components, so that the product is zero. In case one pulls back to
the same component $(X^g)^{\la x,y\ra}$, the obstruction bundle is
equal to that of $\Korb(X^g,Z(g))$, since the respective maps are
given by $e_i:(X^g)^{\la x_1,x_2\ra}\to (X^g)^{x_i}$.
\end{proof}

\begin{cor}
Given $(X,G)$ and $(Y,G)$, $X\times Y$ has a diagonal $G$ action and
$\Korb(I(X\times Y,G),G)=\Korb(I(X,G),G)\widehat{\widehat{\otimes}}
\Korb(I(Y,G),G)$ with the diagonal product structure.
\end{cor}
\begin{proof}
Using the Proposition \ref{inertiacalcprop} above, Proposition
\ref{kunnethprop}
 and the definition of $\widehat{\widehat\otimes}$
\begin{eqnarray}
\Korb(I(X\times Y,G),G)&=&\bigoplus_{g\in G} \Korb((X\times Y)^g
,Z(g))\nn\\
&=& \bigoplus_{g\in G} \Korb(X^g,Z(g))\widehat
\otimes \Korb(Y^g,Z(g))\nn\\
&=&\Korb(I(X,G),G)\widehat{\widehat \otimes} \Korb(I(Y,G),G)\nn\\
\end{eqnarray}
\end{proof}

\begin{cor}
Denote the set of double conjugacy classes of $G\times G$ by $C^2(G)$.
Additively:
\begin{eqnarray}
\label{ksmalleq}
\Ksmall(I(X,G),G)&=&\Korb(I(X,G),G)^G\nn\\
&=&[\bigoplus_{(g,x)\in G\times G, x\in Z(g)} K((X^g)^x)]^{G} \nn\\
&=& \bigoplus_{[g,x]\in C^2(G), x\in Z(g)} K(X^{\la g,x
\ra})^{Z(g,x)}
\end{eqnarray}
and as rings
\begin{equation}
\Ksmall(I(X,G),G)=\bigoplus_{[g]\in C(G)} \Ksmall(X^g,Z(g))
\end{equation}
\end{cor}
\qed

\begin{rmk}
On the other hand we have {\em additively}

\begin{eqnarray}
\label{kfulleq}
\Kfull([X/G])&=&\bigoplus_{[g]\in C(G)}K([X^g/Z(g)])\nn\\
&=&\bigoplus_{[g]\in C(G)}K_{Z(g)}(X^g)\nn\\
&=&\bigoplus_{[g]\in C(G),[x]\in C( Z(g))}K((X^g)^x)^{Z(g,x)}\nn\\
&=&\bigoplus_{[g,x]\in C^2(G)} K(X^{\la g,x
\ra})^{Z(g,x)}
\end{eqnarray}
\end{rmk}

\begin{rmk}
Both versions above are hence additively isomorphic to the sum over
double twisted sectors. In particular, if $G$ is {\em Abelian} then
as vector spaces both versions above are additively given by the
direct sum $\bigoplus_{G\times G}(K(X^{\la g,x\ra}))^G$.
\end{rmk}

\subsection{The Second and Third Appearance of the Drinfel'd double}

Before going on to the twisting it will be instructive to work out
the two theories on the simplest example $[pt/G]$. For both
$\Kfull([pt/G])$ and $\Korb(I(X,G),G)$, we find the Drinfel'd
double, be it in different guises.

\begin{prop}
\label{korbinertia}
 $\Korb(I(pt,G),G)=\subdkg$.
\end{prop}
\begin{proof}
By Proposition \ref{inertiacalcprop}

\begin{equation}
\Korb(I(pt,G),G)=\bigoplus_{g\in G}k[Z(g)]= \bigoplus_{g\in G, z\in
Z(g)} k\eltk{g}{x}
\end{equation}
where we have chosen $\eltk{x}{g}$ for the bi--degree $(g,x)$ part.
Note, all the obstruction bundles vanish, since all the normal
bundles vanish and the multiplication is given by
\begin{equation}
\eltk{g}{x} \eltk{h}{y}=
e_{xy*}(e_x^*(\eltk{g}{x})e_y^*(\eltk{h}{y}))=
\delta_{g,h}\eltk{xy}{g}
\end{equation}
so the multiplication is just that of $k[Z(g)]$.
\end{proof}

\begin{cor}
Since $\Korb(I(pt,G),G)$ is a sum of free $Z(g)$ Frobenius algebras
as needed in Definition \ref{dpralgebrainddef} so we can DPR induce
to obtain  $\dkg$.
\end{cor}
\qed

\begin{cor}\label{moritacor}
As  groupoid algebras the $G$--module  $\Korb(I(pt,G),G)$
is Morita equivalent to $\dkg$.\end{cor}

\begin{proof}
If we consider the $G$ action, we see that it permutes the sectors in
a given conjugacy class. So that the $G$--action on a module
is completely determined
via DPR induction. In the groupoid language, $(I(pt,G),G)$ is
the disjoint union of groupoids $[pt/Z(g)]$ and the $G$--action
adds the morphisms $*_g\stackrel{h}{\to}*_{hgh^{-1}}$  where
$*_g$ denotes the different
objects of the groupoid. This is now Morita equivalent to the loop groupoid
of $[pt/G]$ and hence the result follows.
\end{proof}

See \cite{Wil} for similar considerations.

\begin{thm}
\label{ptdprthm} $\Kfull([pt/G])\iso Rep(\dkg)$.
\end{thm}

\begin{rmk}We were informed by C. Teleman, that a similar formula
at least additively for the case
of $[pt/G]$ can be deduced from the work of Freed-Hopkins-Teleman \cite{FHT1,FHT2}.
\end{rmk}

\begin{nota}
In order to do the calculations, we will use the standard notation
\cite{DVVV,DW,DPR}. Let $\Ag$ be a system of representatives of
conjugacy classes in $C(G)$, which we will consider to be indexed by $A$.
Furthermore let $\alpha$ be an irreducible representation of
$Z(\Ag)$, then we get an irreducible representation $\piAa$ of
$\dkg$ by using DPR induction.
\end{nota}

\begin{proof}[Proof of Theorem \ref{ptdprthm}]
For this we notice that the inertia stack ${\mathfrak
I}[pt/G]=\coprod_{[g]\in C(G)} [pt/Z(g)]$ and hence
$$
\Kfull([pt/G])=\bigoplus_{[g]\in C(G)}
K([pt/Z(g)])=\bigoplus_{[g]\in C(G)} K_{Z(g)}(pt) =
\bigoplus_{[g]\in C(G)} Rep(Z(g))
$$

The product is given by
\begin{multline}
\a_{[\Ag]}*\b_{[\Bg]}=\\\sum_{m_1\in [\Ag],m_2\in [\Bg]}
\frac{|Z(m_1m_2)|}{|G|}
\ind^{Z(\la\bm\ra)}_{Z(\la m_3^{-1}\ra)} \left(\res_{Z(\la
m_1 \ra)}^{Z(\la\bm \ra)}(\a_{m_1})\otimes
\res_{Z(\la m_2\ra)}^{Z(\la \bm \ra)}(\b_{m_2})\right)
\end{multline}

Notice that each $Rep(H)$ has a non--degenerate pairing which is
essentially given by the trace:
$$
\eta(\rho_1,\rho_2):=\frac{1}{|H|}\sum_{h\in
H}tr(\rho_1(h))tr(\rho_2^*(h))
$$
and with this pairing there is an sesquilinar isomorphism of the Frobenius
algebras $(K_G(pt),\chi)$ and $(Rep(G),\eta)$.
What we mean by this is that we can compute the structure constants of the
multiplication for a fixed basis of irreducible representations
 using either metric.

 Furthermore Frobenius
reciprocity holds for a subgroup $H\subset K$
$$
\eta_K(\ind_H^K(\rho_1),\rho_2)=\eta_H(\rho_1,\res_H^K(\rho_2))
$$

So that  we obtain
\begin{eqnarray}
&&\eta(\a_{[\Ag]}*\b_{[\Bg]},\nu_{[\Cg]})\nn\\&=&\hskip-1cm
\sum_{\footnotesize \begin{array}{c} [m_1,m_2,m_3]\\ m_1 \in [\Ag],m_2\in [\Bg],\\
m_3\in [\Cg], \prod m_i=1 \end{array}} \hskip -1cm \eta_{Z(\la
m_3\ra)}(\ind^{Z(\la\bm\ra)}_{Z(\la m_3^{-1}\ra)} \left(\res_{Z(\la
m_1 \ra)}^{Z(\la\bm \ra)}(\a_{m_1})\otimes
\res_{Z(\la m_2\ra)}^{Z(\la \bm \ra)}(\b_{m_2})\right),\nu_{m_3})\nn\\
&=&\hskip-1cm  \sum_{\footnotesize \begin{array}{c}[m_1,m_2,m_3]\\  m_1 \in [\Ag],m_2\in [\Bg],\\
m_3\in [\Cg], \prod m_i=1 \end{array}}  \hskip-1cm \eta_{Z(\la \bm
\ra)}( \res_{Z(\la m_1 \ra)}^{Z(\la\bm \ra)}(\a_{m_1})\otimes
\res_{Z(\la m_2\ra)}^{Z(\la \bm \ra)}(\b_{m_2}),
\res^{Z(\la\bm\ra)}_{Z(\la m_3\ra)}(\nu_{m_3}))\nn\\
&=&\sum_{\footnotesize \begin{array}{c} [m_1,m_2,m_3]\\ m_1 \in [\Ag],m_2\in [\Bg],\\
m_3\in [\Cg], \prod m_i=1\\h\in Z(g_1,g_2) \end{array}}
\frac{1}{|Z(m_2,m_2)|}tr(\a_{m_1}(h))tr(\b_{m_2}(h))tr(\nu_{m_3}^*(h))\nn
\end{eqnarray}
\begin{eqnarray}
&=&\frac{1}{|G|}\sum_{\footnotesize
\begin{array}{c}  m_1 \in [\Ag],m_2\in [\Bg],\\
m_3\in [\Cg], \prod m_i=1\\h\in Z(g_1,g_2) \end{array}}
tr(\a_{m_1}(h))tr(\b_{m_2}(h))tr(\nu^*_{m_3}(h))\nn\\
\end{eqnarray}
for the three--point functions, which agrees with the three point
functions in the case of the Drinfel'd double calculated in
\cite{DPR}. The two point functions then also coincide, since we can
take one representation to be identity, viz.\ the trivial
representation on the identity sector.

\end{proof}

\section{Twisting}
\label{twistpar} In this section, we will be concerned with twisting
of the above structures. This can actually be done on three levels
in two different but equivalent fashions. For the twisting, we can
concern ourselves as above with $(X,G)$ and $(I(X,G),G)$, where we
will consider twisting $\Korb(X,G), \Korb(I(X,G),G)$ and
$\Kfull([X/G])$. The first two are of course isomorphic to the
global orbifold Chow ring or Cohomology ring.

\subsection{Geometric twisting: Gerbe twisting}

In this subsection, we give a geometric interpretation of
the twistings in terms of gerbes.

\begin{assumption}
We will only consider global quotients $(X,G)$ and gerbes
equivariantly pulled back from a point. This means in particular
that we can think of $0,1,2$ gerbes as elements in
$Z^{1,2,3}(G,\uk)$. These gerbes are necessarily flat.
\end{assumption}

\begin{rmk}
It is well known that there is a transgression of an $n$--gerbe on a
stack $\mathfrak X$ to an $n-1$ gerbe on its inertia $\mathfrak{IX}$
\end{rmk}

\subsection{Line bundle twisting}
\label{lintwistsection} Given a line bundle $\L_Y$ on $Y$ there are
basically two ``twists'' we can do. One in K--theory and one in
cohomology, which are as follows. For cohomology, we can consider
cohomology with coefficients in the line bundle $H^*( Y , \L_Y)$ and
in K-theory, we have an endomorphism.
\begin{equation}
K(Y)\stackrel{\sim}{\to} K(Y), \F \mapsto \F\otimes \L_Y
\end{equation}

We will use the notation $K(Y)_{\L}$ to denote the twisted
side.

\begin{rmk}
One way to view this is that the line bundles $\L$ are
gauge degrees of freedom.
\end{rmk}

If we can choose a global section $s$ of $\L$ then we get an
isomorphism
\begin{equation}
\label{sectionisoeq} H^*(Y,k) \to H^*(Y,\L); v\mapsto v \cdot s
\end{equation}

Given  line bundles $\L,\L'$ and $\L''$ on $Y$ and an isomorphism
$\mu:\L\otimes \L'\to \L''$, we get the following multiplicative
maps.
\begin{eqnarray}
H^*(Y,\L)\otimes H^*(Y,\L')&\stackrel{\cup}{\to}&H^*(Y,\L\otimes
\L')\stackrel{\mu_*}{\to} H^*(Y,\L'')\nn\\
K(Y)_{\L}\otimes K(Y)_{\L'}&\to& K(Y)_{\L\otimes \L'} \to K_{\L''}(Y)\nn\\
(\F\otimes \L)\otimes (\F'\otimes \L') &\mapsto&\F\otimes \F'\otimes
(\L\otimes \L')\to\F\otimes \F'\otimes \L''
\end{eqnarray}

\begin{rmk}
If $Y$ has a $G$ action and the line bundles are equivariant line
bundles, then the maps above carry over to the $G$--equivariant
case.
\end{rmk}

\begin{cav}
The equation (\ref{sectionisoeq}) in the $G$--equivariant setting is
only an isomorphism on the level of vector spaces. If the bundle
$\L$ is trivial but the $G$--module structure is given by a
character $\chi$ then the $G$--module structure will be twisted by
$\chi$ upon tensoring with $\L$.
\end{cav}

\subsection{0-Gerbe twisting: Ramond twist}
By definition a $0$--gerbe is nothing but a line bundle on the stack
and if we are dealing with a global quotient $(X,G)$, using the
assumption above, we get a trivial line bundle $\mathscr L$ on $X$,
which is equivariant, but not necessarily equivariantly trivial.

If we fix a trivialization of the line bundle, viz.\ choose a global
section $v$. This induces an isomorphism

\begin{equation}
\mu:\mathscr{L}\otimes \mathscr{L}\to \mathscr{L};\quad  v\otimes
v\mapsto v
\end{equation}

The equivariance of this line bundle is expressed by isomorphisms

\begin{equation}
g^*(\L)\iso \L; \quad  v\mapsto \chi(g)v; \quad \chi\in
Z^1(G,\uk)=Hom(G,\uk)
\end{equation}

In terms of the twisting using $\mu$, we can twist as described in
the paragraph above. In this case, the $G$--action will be twisted
by the character $\chi$ as will be the metric. This will ``destroy''
the properties of a pure $G$--Frobenius algebra (for instance axiom
T will cease to hold), but we will almost end up with a
$G$--Frobenius algebra which is twisted by the character $\chi$.
This will indeed be the case, if we had started out with a
$\chi^{-1}$ twisted Ramond model \cite{wisc,orb}. In the current
A--model setting, we will always have invariant metrics and strict
self--invariance (axiom T). This type of twist is, however, very
important in the B--model setting as it is not guaranteed that the
objects have invariant pairings and self--invariance \cite{wisc,orb,
singorb}. Hence we can view the 0--gerbe twisting as a twisting to
the Ramond model and hence as a spectral flow \cite{DVV,wisc,orb,
singorb}.

\subsection{1-Gerbe twisting: discrete torsion}
\label{onegerbesection}
This twisting has been investigated the most and goes under the name
of discrete torsion. We shall disentangle the definitions so as to
show that the resulting algebraic structure is that of \cite{disc}.
This exposition owes a lot to \cite{thaddeus} and \cite{hitch}.

A 1--Gerbe $\G$ on $(X,G)$ which is equivariantly pulled back from a
point is given by fixing the (a) isomorphism $\L_g:g^*(\G)\isoto\G$
which are in turn given by line bundles $\L_g$ and (b) isomorphisms
$\psi(g,h):\L_g\otimes \L_h\to \L_{gh}$ which are associative.

Notice since the gerbe is trivial on $X$, so are the line bundles.
In order to go on, we also choose sections $s_g$ of $\L_g$. Then in
this basis the morphisms $\psi(g,h)$ are given by their matrix entry
$\a(g,h)\in Z^2(G,\uk)$. Notice that a different choice of sections
changes $\a$ by a co-boundary.

\begin{rmk}
Notice that the line bundles $\L_g|_{X^g}$ are actually $Z(g)$
equivariant line bundles. Furthermore fixing the sections $s_g$ we
see that the isomorphisms are given by the characters
$\eps_g(h)=\a(g,h)/\a(h,g)$ which are the famous discrete torsion
co--cycles (see the e.g.\ \cite{disc} for a full list of
references). Furthermore, $\eps(g,h):=\eps_g(h)$ is even a
bi--character when restricted to commuting elements (see e.g.\
\cite{disc}). This means that as $Z(m_1)\cap Z(m_2)$ modules
$\L_{m_1}\otimes \L_{m_2}|_{X^{\bm}}\iso \L_{m_1m_2}|_{X^{\bm}}$.
\end{rmk}

\subsubsection{Cohomology}

We can now set ${\mathscr H}^{\G}(X,G):=\bigoplus
H^*(X^g,\L_g|_{X^g})$. For the multiplication, we can use the
standard push-pull mechanism in a slightly modified version: for
$v_{m_i}\in H^*(X^{m_i},\L_{m_i}|X^{m_i})$
\begin{equation}
v_{m_1}*_{\G}v_{m_2}:=\check
e_{m_3*}(\psi_*(m_1,m_2)_{X^{\bm}}[e^*_1(v_{m_1})e_2^*(v_{m_2})]\e(\R(\bm)))
\end{equation}
Notice that the result indeed lies in
$H^*(X^{m_1m_2},\L_{m_1m_2}|_{X^{m_1m_2}})$ due to the projection
formula.

Given the section $s_g$ we get isomorphism of the
$\lambda_g:H^*(X^g,\L_g|_{X^g})\isoto H^*(X^g)$ additively  and this
induces a new twisted multiplication on $A:=\bigoplus H^*(X^g)$ via

\begin{eqnarray}
&&v_{m_1}*_{\a}v_{m_2}\nn\\&:=&\lambda_{m_1m_2}^{-1}\circ \check
e_{m_3*}(\psi_*(m_1,m_2)|_{X^{\bm}}[e^*_1(\lambda_{m_1}(v_{m_1}))
e_2^*(\lambda_{m_1}(v_{m_2}))]\e(\R(\bm))\nn\\
&=&\a(m_1,m_2)\, v_{m_1}*v_{m_2}
\end{eqnarray}
That is we realize the algebraic twist of \cite{disc} and \S \ref{algtwistsec}
above.

Of course we could have alternatively discussed the Chow ring $A^*$ in
the same way.

\subsubsection{K-theory I: twisted multiplication}
In the case of $K$-theory using the standard formalism, we will
obtain morphisms

\begin{equation}
K(X^{m_1})_{\L_{m_1}|_{X^{m_1}}}\otimes
K(X^{m_2})_{\L_{m_2}|_{X^{m_2}}}\rightarrow
K(X^{m_1m_2})_{\L_{m_1m_2}|_{X^{m_1m_2}}}
\end{equation}

Considering the direct sum of twisted $K$--theories
\begin{equation}
\Kglobal^{\G}(X,G):=\bigoplus_{m \in G} K(X^m)\otimes \L_m
\end{equation}
we hence obtain a multiplication using the push--pull
formalism of equation (\ref{proddefeq}) analogously to the above.

And by choosing sections, we again get a twisted version of the
multiplication
\begin{equation}
\F_{m_1}*_{\a}\F_{m_2}=\a(m_1,m_2)\, \F_{m_1}*\F_{m_2}
\end{equation}
where a different choice of sections results in a change of $\a$ by
a co--boundary.

\begin{rmk} There are several aspects, though not all,
of the considerations above which have been previously discussed
and also there have been related discussions which we would like to address briefly:
\begin{itemize}
\item It was shown in \cite{AR1} that the {\em additive} $\a$
twisted $K$--theory of $(X,G)$ as defined via projective
representations is given by $K^{\a}\iso \bigoplus_{[g]}(
K(X^g)\otimes \L_g)^{Z(g)}$ where $\L_g$ was considered as a $G$
module via the discrete torsion co--cycle
$\eps(g,h)=\a(g,h)/\a(h,g)$ for $[g,h]=e$. There is no obvious
multiplicative structure on this space as remarked in \cite{AR1},
but the formalism above does give it a multiplicative structure.

\item We would also like to note that in \cite{LU} an additive theory
for a gerbe twist was constructed and it was shown
that in the case of a global
quotient with a gerbe pulled back from a point the gerbe twisted $K$--theory
and the Adem--Ruan twisted theory as cited above coincide.

\item Our geometric twisting above coincides
with the algebraic twisting of
GFAs considered in \cite{disc} and \cite{JKK2} ---  see \S\ref{algtwistsec} below.
Hence the formula above and the Chen character of \cite{JKK2} answer the
question of Thaddeus \cite{thaddeus} about the relation of the two types of possible twists
by line bundles in Cohomology vs.\ K--theory.

\end{itemize}
\end{rmk}

\subsubsection{K--theory II: twisted K-theory}
Another standard thing to do with a flat gerbe, that is a 2--cocycle
$\theta\in H^2(Y,\uk)$ is to regard the twisted $K$--theory
$K^{\theta}(Y)$. In our case of a global orbifold, given $\a$ as
above we will study the twisted equivariant $K$--theory
$K_G^{\a}(X)$ which by definition is the twisted $K$--theory of the
stack $K^{\a}([X/G])$.

In this interpretation, one cannot see any type of multiplication.
It is basically the same problem as in the case of a 0--gerbe. The
natural product goes from $K^{\a}(Y)\otimes K^{\b}(Y)\to K^{\a
\beta}(Y)$. We will get back to this in the 2--gerbe twisting.

\subsubsection{Twisted group ring}
It is again  useful to look at the details in the case of $(pt,G)$.
Here $K^{\a}([pt/G])=Rep^{\a}(G)$ that is the ring of projective
representation with cocycle $\a$.

On the other hand the global orbifold K--theory with an $\a$ twist
$\Ksmall^{\a}(pt,G)=k^{\a}[G]$ and the $G$ invariants by the
conjugation action are isomorphic to $Rep^{\a}(G)$ \cite{Kar}.

Here the multiplication is the one in $k^{\a}[G]$ which is just the
one of $k[G]$ twisted by $\a$.

\subsection{2-Gerbe twisting}
Finally, we wish to discuss twisting by a gerbe of the type $\b\in
Z^3(G, \uk)$. This type of gerbe is also the one we used to twist
the Drinfel'd double and indeed there is a connection.

We can transgress the equivariant 2--gerbe to an equivariant
1--gerbe $\G$ on $\mathfrak IX$ and actually even to a 1--gerbe over
$(I(X,G),G)$. Here the gerbe is characterized by a set of line
bundles, which provide the isomorphisms $\L_{g,x}:
x^*(\G|_{X^g})\isoto \G|_{X^{x^{-1}gx}}$ together with associativity
isomorphisms $\theta_g(x,y):\L_{g,x}\otimes \L_{h,y}\to \L_{g,xy}$
if $g=x^{-1}gx$.

The condition of coming from a 2-gerbe expresses itself in a
constraint on the $\theta_g$. In particular it means (see e.g.\
\cite{Wil}) the $\theta_g$ are given by equation (\ref{thetaeq}).

\subsubsection{2--Gerbe twisted K--theory I: twisting on $\Korb((I(X,G),G))$}
Now we are in a situation in which we can twist.

First of all there is a na{\"i}ve twisting on
$\Korb((I(X,G),G)$ by the various $\theta_g$ transgressed from $\beta$;
see \S \ref{twistingtwopar} below where we give a more detailed description
of this type of twist.
In the case of $(pt,G)$ with $G$ Abelian this yields a geometric
incarnation of $\dbkg$. In the general group case, we get a Morita equivalent
subalgebra just as in the untwisted case.

\subsubsection{2--Gerbe twisted K--theory II: twisting on $\Kfull([X/G])$}
More importantly, however, there is a twisting for the full K--theory.

\begin{df}
Given $\b \in Z^3(G,(\uk)$
we define the twisted full $K$--theory $\Kfull^{\beta}([X/G])$
 using the co--product and the obstruction:
that is the multiplication which is induced by:
\begin{equation}
\label{betamulteq}
\cf_{g} \cdot \cf_{h} :=   e_{3*}(e_{1}^*(\cf_{g}) \otimes^{\gamma}
e_{2}^*(\cf_{h}) \otimes Obs_K(g,h))
\end{equation}
see \cite{JKK2} for details on how this global formula relates to the
inertia stack setting.

Here we use the co--product in $\dbkg$ which is given by $\gamma$ defined above
by equation (\ref{gammadefeq}) to define the action of $Z(g,h)$ on the
tensored bundle. This means that if for $x\in Z(g,h)$
 $\phi_x:x^*(\cf_g)\to \cf_{g}$ and $\psi_x:x^*(\cf_h)\to \cf_{h}$
  are the isomorphisms given by the equivariant data, then
the isomorphism of $x^*(e_{1}^*(\cf_{g}) \otimes^{\gamma}
e_{2}^*(\cf_{h}))\iso e_{1}^*(\cf_{g}) \otimes^{\gamma}
e_{2}^*(\cf_{h})$ is chosen to be $\gamma_x(g,h) \phi_x|_{X^{g,h}}
\otimes \psi_x|_{X^{g,h}}$ where $\gamma$ is defined by  equation
(\ref{gammadefeq}).

\end{df}

\begin{rmk}
For an interesting, different and independent approach we refer the
reader to \cite{AR2}. Here the authors consider a twist which is on
the full K-theory of the inertia stack $\Kfull(\ccx)$ and does not
seem to use a co--product structure. The latter is key to the
braided associativity.
\end{rmk}

\subsection{Case of a point}
Restricting to a point, we obtain the
analog of Theorem \ref{ptdprthm}:

\begin{thm} \label{dbkgrepthm}
 $\Kfull^{\beta}([X/G])^{\beta}\iso \dbkg$
\end{thm}

\begin{proof}
Analogous to Theorem \ref{ptdprthm} using the calculations of
 \cite{DPR,DW,DVVV}.
\end{proof}
\begin{rmk}
Here we see that we essentially get the \cite{DW} realization of the
2-d calculation of \cite{DVVV}, which is astonishing and inspiring.
Using this insight, we can also understand why the twist already
works on the level of the global quotient stack itself. The point is
that applying the full stringy $K$--theory functor already entails
moving to the inertial stack. This can be interpreted as moving to
the loop space and hence evaluating the correlation functions on
$\sigma\times S^1$, viz.\ the procedure described in \cite{DW}. This
explains why the 1+1 dimensional theory has the flavor of a 2+1
dimensional theory.
\end{rmk}

\begin{rmk}
This theorem is mathematically
astonishing in the sense that the resulting structure is
neither commutative nor associative in general.
{\em We will get an essentially non--associative
algebra unless $\beta\equiv 1$.}
But it is of course
associative and commutative in the sense
of braided monoidal categories.
We hope that we have motivated the
appearance of braided monoidal categories
already through the definition of
Frobenius traces and objects. Moreover,
if one reads for instance Moore and Seiberg's work on
classical and quantum
field theory one sees that the fusion ring is actually
not expected to be
associative and commutative. However, the fusion and braiding operators satisfy
pentagon and hexagon relations.
Only the dimensions of the intertwiners lead
to such an algebra on the nose. In case of the objects themselves one
should actually expect that one has to go to the braided picture.
\end{rmk}

\subsubsection{Verlinde algebra}
We can get an associative algebra by introducing a basis
of irreducible representations $V_i, i\in I$
and using
the dimensions of the intertwiners as the structure coefficients.
That is if $V_i\otimes V_j=\bigoplus_k V_{ij}^k \otimes V_k$,
where  $V_{ij}^k$ is the space of intertwiners or multiplicity,
set $c_{ij}^k= dim(V_{ij}^k)$. Then the Verlinde ring is just
$k[v_i,i\in I]$ where the $v_i$ are now formal variables with the
multiplication $v_iv_j=\sum_k c_{ij}^k v_k$.

\subsection{Algebraic twisting}
\label{algtwistsec}
In this section, we give a purely algebraic version of the twistings.
This allows us among other things
 to connect the 1--Gerbe twistings to
the discrete torsion twistings used in \cite{disc,JKK2}.

\subsubsection{Algebraic Twisting I: Discrete Torsion}
We briefly recall the twisting by discrete torsion in the
$G$--Frobenius algebra case.

In \cite{disc} we defined the twisting of $G$--Frobenius algebras
via
\begin{equation}
A\leadsto A^{\a}:=A\widehat\otimes k^{\a}[G]
\end{equation}
This provides an action of the group $Z^2(G,\uk)$ on the set of
GFAs. Notice that two twists $A^{\a}$ and $A^{\b}$ are
isomorphic if and only if $[\a]=[\b]\in H^2(G,\uk)$. It is clear that
this extends to $G$-Frobenius algebra objects.

\begin{prop}
\label{twistequalprop} The algebraic twist and the geometric twist
coincide, that is  for $\a \in Z^2(G,\uk)$
\begin{equation}
(\Forb(X,G))^{\a}=\Forb^{\a}(X,G)
\end{equation}
\end{prop}
\begin{proof}
Straightforward from the definition and paragraph \S\ref{onegerbesection}.
\end{proof}

\subsubsection{Algebraic twisting II: Twisting on $I(X,G)$ and the second appearance of the twisted Drinfel'd double}
\label{twistingtwopar}

Notice that by Proposition \ref{inertiacalcprop}
$K_{global}(I(G,X),G)$ splits as a direct sum of rings indexed by
$g\in G$, each of which is a $Z(g)$--Frobenius algebra. If $G$ is
Abelian then all the $Z(g)=G$. It is hence possible to twist each
$G$--Frobenius algebra separately by a discrete torsion $\theta_g\in
Z^2(G,\uk)$. In the non--Abelian case, the twists can not be chosen
arbitrarily, since they have to be compatible with the $G$ action
that acts by double conjugation. This means that one has the free
choice of a twist for each conjugacy class $[g]$, that is co-cycles
$\theta_g\in Z^2(G,\uk)$, such that
$\theta_g(h,k)=\theta_{xgx^{-1}}(xhx^{-1},xkx^{-1})$ for all $x\in
X$.

In this situation, we can also ask that the $\theta_g$ be even more coherent,
that is that they stem from a $\beta\in Z^{3}(G,\uk))$. In this case
we  basically obtain an identification
of $\Korb(I(pt,G),G)$ with the Drinfel'd double.

\begin{dfprop} For $\beta\in Z^3(G,\uk)$ 
\begin{eqnarray}
\Korb^{\beta}(I(X,G),G)&:=&\bigoplus_{g\in G}\Korb^{\theta_g}(X^g,Z(G))\\
&=&\bigoplus_{g\in G}\Korb(X^g,Z(G))\hotimes
k^{\theta_g}[Z(g)]\\
\label{propparteq}&=&\Korb(I(X,G),G)\hhotimes \dbkg\\
\end{eqnarray}
\end{dfprop}
\begin{proof}
The proposition part is the equation (\ref{propparteq}). In view of
Proposition \ref{twistequalprop} this follows from the fact that the
components $(g,x)$ of $\Korb$ are only non--empty if $x\in Z(g)$.
The restriction to the corresponding subspace of $\dbkg$ is given by
$\bigoplus_g k^{\theta_g}[Z(G)]$.

\end{proof}

\begin{cor}
There is an action of $Z^3(G,\uk)$ on $\Korb(I(X,G),G)$ obtained
by tensoring with $\hhotimes \dbkg$.
\end{cor}

\begin{proof}
Directly from the above and Lemma \ref{hhotimeslem}.
\end{proof}

We can thus twist $\Korb(I(X,G),G)$
via the procedure above and hence have a completely analogous
story to the twists of $\Korb(X,G)$ by discrete torsion analyzed
in \cite{disc}, but
now one gerbe level higher.

If $\Korb(I(X,G),G)$ is free in the sense
that all the $Z(g)$--Frobenius algebras are free, we can further
 DPR--algebra induce as discussed in \S\ref{dprindsection}.

\begin{thm}
We have the following identifications:
\begin{equation}
(\Korb^{\beta}(I(pt,G),G))=\subdbkg
\end{equation}
 and
\begin{equation}
\dprind(\Korb^{\beta}(I(pt,G),G))=\dbkg
\end{equation}
\end{thm}

\begin{proof}
Straightforward computation.
\end{proof}

\subsubsection{Algebraic twisting III}
In contrast to the previous twistings, the full orbifold $K$--theory
twisting cannot just be reduced to an algebraic twisting. This can
only be done additively in general. In the trival $G$--action case however, the
twists by $\beta\in Z^3(G,\uk)$ again has a purely algebraic
description.

\begin{prop}
\label{gltwprop}
Given a global quotient stack $\cx=[X/G]$
and a class $\beta \in Z^3(G,\uk)$,
we have {\em additively}
\begin{eqnarray}
\Kfull^{\beta}(\cx)&:=&\bigoplus_{[g]} K^{\theta_g}[X^g/Z(g)]\\
&=&\bigoplus_{g}\left( \Korb((X^g,Z(g)))\hotimes k[Z(g)]\right)^{Z^g}\\\end{eqnarray}
but the multiplication
is the one defined by equation \ref{betamulteq}.
\end{prop}

\begin{proof}
This follows from the fact that additively $K_H(Y)\iso
\Kglobal((Y,H))^H$.
\end{proof}

\subsubsection{Trivial action case}
In the case of a trivial  $G$--action, the multiplication becomes
particularly transparent.

\begin{thm}
\label{trivthm}
Let $\cx=[X/G]$ where $X$ has a trivial $G$ action then:
\begin{eqnarray}
\Kfull^{\beta}(\cx)&\iso&
 \bigoplus_{[g]} K(X)\otimes Rep^{\theta_g}(Z(g)) \\
&\iso&K(X)\otimes Rep((\dbkg))
  \end{eqnarray}
where in the last line the algebra structure is the tensor product and
 in the second line we have the following multiplication:

\begin{equation}
\cf_g\otimes \rho * \cf_h \otimes \rho':=
\cf_g*\cf_h\otimes \rho*\rho'
\end{equation}
where $\cf_g*\cf_h=\cf_g\otimes \cf_h\in K(X)$
 and $\rho*\rho'$
is induced by
$\res^{\dbkg}_{Z(gh)}( \dprind(\rho)\otimes \dprind(\rho'))$
using the braided structure of $\dbkg$ and Theorem \ref{DPRthm}.
\end{thm}

\begin{proof}
First we calculate using that for a trivial action all $X^g=X$:
\begin{eqnarray}
\Kfull^{\beta}(\cx)&=&\bigoplus_{[g]} K^{\theta_g}[X^g/Z(g)]\\
&\iso&\bigoplus_{[g]} K_{Z(g)}^{\theta_g} (X)\\
&\iso&\bigoplus_{[g]} K (X)\otimes Rep^{\theta_g}(Z(g)) \\
&\iso& K (X)\otimes \bigoplus_{[g]} Rep^{\theta_g}(Z(g)) \\
&\iso&K(X)\otimes Rep((\dbkg))
 \end{eqnarray}
Where the second line is by Grothendieck, the
third line follows from e.g.\ Lemma 7.3 of \cite{AR1}
and the forth line uses Theorem \ref{DPRthm}.

Now for the multiplicative structure, we notice that since
the inclusions $e_i$ are all the identity, on the
factors of $K(X)$ the multiplication boils down to the tensor product,
whereas the product on the representations rings goes through
the induction process and uses the co--cycle $\gamma$. This
of course is nothing but the description in terms of $\dbkg$.
\end{proof}

\subsection{Alternative description using modules}

The Theorem \ref{trivthm} above can nicely be seen in the module
language.

\subsubsection{Equivariant $K$--theory in the module language}
We first recall the setup of $G$--equivariant $K$--theory in terms
of modules, see \cite{Atiyah,AS}: $K_G(X)\iso B_{\rm proj., fin.\,
gen.}-mod$ where $B$ is $C^{\infty}(X)\rtimes G$ with the
multiplication $(a,g)\cdot (a',g')=(ag(a'), gg')$ and the modules
are projective finitely generated.

In order to twist with a 1--gerbe $\a\in Z^2(G,\uk)$ following Atiyah-Segal,
we give a new multiplication on $B$ via

\begin{equation}
\label{alphatwisteq}
(a,g)\cdot (a',g')=(ag(a'),\a(g,g') gg')
\end{equation}
We call the resulting ring $B^{\a}$. Now the twisted $K$--theory is
given by the projective finitely generated $B^{\a}$--modules.

The na\"ive tensor structure which sends the $\a$ twisted $K$--theory times
the $\beta$ twisted $K$--theory to the $\a\b$ twisted $K$--theory
 uses the $A$ module structure induced
by the multiplication map $A\otimes A\to A$
and the co--product $\Delta:k[G]\to k[G]\otimes k[G]$ given
by $\Delta(g)=g\otimes g$.

\begin{rmk}
In the algebraic category, we
can use $\mathcal O_X$ instead of $C^{\infty}(X)$.
\end{rmk}

\subsubsection{Remarks on the 2--gerbe twisted case}

For the 2-gerbe $\b$ twisted $K$--theory, we can describe the
$K$--theory additively as follows. Let $A_g=C^{\infty}(X^g)$, let
the $\theta_g$ be defined via equation \ref{thetaeq} and we define
$B_g^{\theta_g}=A_g\rtimes Z(g)$ with the multiplication as in
equation \ref{alphatwisteq}. Set  $B=\bigoplus_{[g]} B_g^{\theta_g}$
then additively $\Kfull^{\beta}([X/G])\iso  B_{\rm proj., fin.\,
gen.}-mod$

To describe the multiplicative structure we would have to define a
co--product on $B$ which incorporates  the $G$--grading, the
obstruction and the twisting. This should also be possible for a
general stack or groupoid and a 2--gerbe. The full analysis is
beyond the scope of the present considerations, but we plan to
return to this in the future.

In the special case of a trivial $G$
action the construction has can be made fully explicit.

\subsubsection{Trivial $G$--action}
In the trivial $G$ action case, like  the case $[pt/G]$, there are
no obstructions and we can give a full description of the theory in
module terms: Let $A=C^{\infty}(X)$ {\em and} assume $X$ has a
trivial $G$ action then $C^{\infty}(I(X,G))=\bigoplus_{g\in G} A$.
Now although the $G$ action on $X$ is trivial, it is not trivial on
$I(X,G)$, since it permutes the components.

It is easy to check that in the trivial $G$ action case the algebra is
\begin{equation}
B=\bigoplus_{[g]\in G} C^{\infty}(X)\rtimes k[Z(g)] \morita A\otimes
\dkg
\end{equation}
where the product structure on the $K$--theory is given via the
co--multiplication of $\dkg$.

Similarly, twisting with $\b$ we obtain

\begin{prop}
$$
B^{\beta}\morita A\otimes \dbkg
$$
and the product structure on the $K$--theory of projective finitely
generated $B^{\b}$--modules is given via the co--multiplication of
$\dbkg$.
\end{prop}

By considering the braided category projective finitely generated
$B^{\beta}$ modules we hence obtain  a generalization of the Theorem
of $[pt/G]$ to the case of a trivial $G$--action  which is analogous
to Theorem \ref{trivthm}.

\end{document}